\documentclass[aop,seceqn,citesort,MSNbibl,dvips]{arximspdf}

\doi{10.1214/10-AOP552}
\volume{39}
\issue{1}
\pubyear{2011}
\firstpage{369}
\lastpage{405}

\makeatletter

\newtheorem{Theorem}{Theorem}[section]
\newproclaim{Definition}[Theorem]{Definition}
\newtheorem{Lemma}[Theorem]{Lemma}
\newproclaim{Remark}[Theorem]{Remark}
\newtheorem{Hypothesis}{Hypothesis}

\newcommand{\nat}{\mathbb N}
\newcommand{\E}{\mathbb E}
\newcommand{\Pro}{\mathbb P}
\newcommand{\reals}{\mathbb R}

\makeatother

\begin{document}
\begin{frontmatter}

\title{Fast transport asymptotics for stochastic RDEs with boundary noise}
\runtitle{Fast transport asymptotics for stochastic RDEs}

\begin{aug}
\author[B]{\fnms{Sandra} \snm{Cerrai}\corref{}\thanksref{t1}\ead[label=e1]{cerrai@math.umd.edu}} and
\author[B]{\fnms{Mark} \snm{Freidlin}\thanksref{t2}}
\runauthor{S. Cerrai and M. Freidlin}
\affiliation{University of Maryland}
\address[B]{Department of Mathematics\\
University of Maryland\\
College Park\\
Maryland, 20742\\
USA\\
\printead{e1}}
\end{aug}

\thankstext{t1}{Supported in part by NSF Grant DMS-09-07295.}
\thankstext{t2}{Supported in part by NSF Grant DMS-08-03287.}

\received{\smonth{7} \syear{2008}}
\revised{\smonth{6} \syear{2009}}

%
\begin{abstract}
We consider a class of stochastic reaction-diffusion equations also
having a stochastic perturbation on the boundary and we show that when
the diffusion rate is much larger than the rate of reaction, it is
possible to replace the SPDE by a suitable one-dimensional
stochastic differential equation. This replacement is possible under
the assumption of spectral gap for the diffusion and is a result of
\textit{averaging} in the fast spatial transport. We also study the
fluctuations around the averaged motion.
\end{abstract}

%

\end{frontmatter}

\section{Introduction}
\label{sec1}

In classical chemical kinetics, the evolution of concentrations of
various components in a reaction is described by ordinary differential
equations. Such a description turns out to be unsatisfactory in a
number of
applications, especially in biology (see \cite{volke}).

There are several ways to construct
a more adequate mathematical model.
If the reaction is fast enough,
one should take into account that the concentration is not constant
in space in the volume where the reaction takes place. Then, the change of
concentration due to spatial transport
should be taken into account and the system of ODEs should be replaced
by a
system of PDEs of reaction-diffusion type.
In some cases, one should also take into account random changes in time
of the rates
of reaction. Then, the ODE is replaced by a stochastic differential equation.
If the rates change randomly not just in time but also in space, then
evolution of concentrations can be described by a system of SPDEs.

On the other hand, the rates of chemical reactions in the system and the
diffusion coefficients may, and as a rule do, have different orders. Some
of them are much smaller than others and this allows one to apply
various versions of
the averaging principle and other asymptotic methods, thereby
eventually obtaining a relatively simple description of the system.

In this paper, we study the case where the diffusion rate is much larger
than the rate of reaction and we show that in this case, it is possible
to replace SPDEs of reaction-diffusion type by suitable SDEs. Such
an approximation is valid, in particular, if the reaction occurs only
on the
boundary of the domain (this means that the nonlinearity is included in
the boundary conditions). This replacement is a result of
\textit{averaging} in
the fast spatial transport. We would like to stress that our approach
allows us also to calculate the main
terms of deviations of the solution of the original problem from the
simplified model. Notice, moreover, that the case where the diffusion
coefficients and some of the
reaction rates are large compared with other rates can be considered in a
similar way.

More precisely, we are dealing with the following class of equations:
%
%
\begin{equation}
\label{eq1}
\cases{
\dfrac{\partial u_\varepsilon}{\partial t}(t,x)=\dfrac{1}{\varepsilon
}\mathcal{A} u_\varepsilon(t,x)+f(t,x,u_\varepsilon
(t,x))\cr
\hphantom{\,\dfrac{\partial u_\varepsilon}{\partial t}(t,x)=}
{}+g(t,x,u_\varepsilon(t,x)) \,\dfrac{\partial w^{Q}}{\partial
t}(t,x), &\quad $t\geq0, x \in D$,\cr
\dfrac1\varepsilon\dfrac{\partial u_\varepsilon}{\partial\nu
}(t,x)=\sigma(t,x)\dfrac{\partial w^B}{\partial t}(t,x), &\quad $t\geq0,
x \in\partial D$,\vspace*{2pt}\cr
u_\varepsilon(0,x)=u_0(x), &\quad $x \in D$,}
\end{equation}
for some $0<\varepsilon\ll 1$.
These are reaction-diffusion equations perturbed by a noise of
multiplicative type, where the diffusion term $\mathcal{A}$ is
multiplied by a large parameter $\varepsilon^{-1}$ and a noisy
perturbation is also acting on the boundary of the domain~$D$.

Here, $D$ is a bounded open subset of $\reals^d$, with $d\geq1$,
having a regular boundary (for more details, see Section \ref{sec2})
and, in the case $d=1$, we take $D=[a,b]$. $\mathcal{A}$ is a
uniformly elliptic second order operator and $\partial/\partial\nu$
is the corresponding conormal derivative. This is why the same constant
$\varepsilon^{-1}$, which is in front of the operator $\mathcal{A}$, is
also present in front of the conormal derivative $\partial/\partial
\nu$. In what follows, we shall denote by $A$ the realization in
$L^2(D)$ of the differential operator $\mathcal{A}$, endowed with the
conormal boundary condition.

The coefficients $f,g\dvtx[0,\infty)\times D\times\reals\to\reals$ are
assumed to be measurable and satisfy a Lipschitz condition with respect
to the third variable, uniformly with respect to the first two
variables, and the mapping $\sigma\dvtx[0,\infty)\times\partial D\to
\reals
$ is bounded with respect to the space variable.

The noisy perturbations are given by two independent cylindrical Wiener
processes, $w^Q$ and $w^B$, defined on the same stochastic basis
$(\Omega, \mathcal{F}, \mathcal{F}_t, \Pro)$, which take values on
$L^2(D)$ and $L^2(\partial D)$, respectively, and have covariance
operators $Q \in\mathcal{L}^+(L^2(D))$ and $B \in\mathcal
{L}^+(L^2(\partial D))$,
respectively.\setcounter{footnote}{1}\footnote{Here, and in what
follows, given any Banach space $E$, we denote by $\mathcal{L}(E)$ the
Banach space of bounded linear operators on $E$ and by $\mathcal
{L}^+(E)$ the subspace of nonnegative and symmetric operators.} In
space dimension $d=1$, we can take $Q$ equal to the identity operator
so that we can deal with space--time white noise. Moreover, as $L^2(\{
a,b\})=\reals^2$, in space dimension $d=1$, we do not assume any
condition on $B$.

Stochastic partial differential equations with a noisy term also acting
on the boundary have been studied by several authors; see, for example,
da Prato and Zabczyk \cite{dapratozab}, Freidlin and Wentzell
\cite{freidlinwentzell} and Sowers \cite{sowers}. The last two mentioned
papers also deal with some limiting results with respect to small
parameters appearing in front of the noise.
However, the limiting results which we are studying in the present
paper seem to be completely new and we are not aware of any previous
results dealing with the same sort of multiscaling problem, even in the
simpler case of homogeneous boundary conditions (i.e., $\sigma=0$).

As mentioned above, our interest is in studying the limiting behavior
of the solution $u_\varepsilon$ of problem (\ref{eq1}) as the parameter
$\varepsilon$ goes to zero, under the assumption that the diffusion
$X_t$ associated with the operator $\mathcal{A}$, endowed with the
conormal boundary condition [this corresponds to a diffusion $X_t$ on
some probability space $(\hat{\Omega}, \hat{\mathcal{F}}, \hat
{\mathcal{F}}_t, \hat{\Pro})$ which reflect on the boundary of
$D$], admits a unique invariant measure $\mu$ and a spectral gap
occurs. That is, for any $h \in L^2(D,\mu)$,
\[
\int_D \biggl| \hat{\mathbb{E}}^xh(X_t)-\int_Dh(y) \mu(dy)\biggr|^2
\mu(dx)\leq c e^{-2\gamma t}\int_D |h(y)|^2 \mu(dy)
\]
for some constant $\gamma>0$.
This can be expressed in terms of the semigroup $e^{tA}$ associated
with the diffusion $X_t$, by saying that
%
%
\begin{equation}
\label{dubbio}
\biggl|e^{tA} h-\int_D h(x) \mu(dx)\biggr|_{L^2(D,\mu)}\leq c
e^{-\gamma t}|h|_{L^2(D,\mu)}.
\end{equation}
Moreover, as shown in Remark \ref{rem21}, the space $L^2(D)$ is
continuously embedded into $L^2(D,\mu)$.

Our aim is to prove that equation (\ref{eq1}) can be replaced by a
suitable one-dimensional stochastic differential equation, whose
coefficients are obtained by averaging the coefficients and the noises
in (\ref{eq1}) with respect to the invariant measure $\mu$. More
precisely, for any $h \in L^2(D,\mu)$, we define
\[
\hat{F}(t,h)=\int_Df(t,x,h(x)) \mu(dx),\qquad t\geq0,
\]
and for any $h \in L^2(D,\mu)$, $z \in L^2(D)$ and $k \in
L^2(\partial D)$, we define
\[
\hat{G}(t,h) z=\int_D g(t,x,h(x))z(x) \mu(dx),\qquad t\geq0,
\]
and
\[
\hat{\Sigma}(t)k=\delta_0\int_DN_{\delta_0}[\sigma(t,\cdot)k](x)
\mu(dx),\qquad t\geq0,
\]
where $N_{\delta_0}$ is the Neumann map associated with $\mathcal{A}$ and
$\delta_0$ is a suitable constant (see Section \ref{sec2},
\cite{lasiecka} and \cite{lionsmagenes} for definitions). We prove that for
any $t\geq0$, the mappings $\hat{F}(t,\cdot)\dvtx L^2(D,\mu)\to\reals$
and $\hat{G}(t,\cdot)\dvtx L^2(D,\mu)\to L^2(D)$ are both well defined
and Lipschitz continuous, and $\hat{\Sigma}(t) \in L^2(\partial
D)$, so that the stochastic ordinary differential equation
%
%
\begin{equation}
\label{eq1bis}
\cases{
dv(t)=\hat{F}(t,v(t)) \,dt+\hat{G}(t,v(t))\,dw^Q(t)+
\hat{\Sigma}(t)\,dw^B(t),\vspace*{2pt}\cr
\displaystyle v(0)=\int_D u_0(x) \mu(dx),}
\end{equation}
admits, for any $T>0$ and $p\geq1$, a unique strong solution $u \in
L^p(\Omega;C([0,T]))$ which is adapted to the filtration of the noises
$w^Q$ and $w^B$. Notice that (\ref{eq1bis}) is a one-dimensional
stochastic equation, in the sense that the space variables have
disappeared. In Section \ref{sec4}, we show that it can be rewritten as
\[
dv(t)=\hat{F}(t,v(t)) \,dt+\Phi(t,v(t)) \,d\beta_t,
\]
where $\beta_t$ is a standard Brownian motion and the diffusion
coefficient $\Phi$
is explicitly given in terms of $Q$, $G$, $B$ and $\Sigma$.

When we say that equation (\ref{eq1}) can be replaced by
(\ref{eq1bis}), we mean that the solution $u_\varepsilon$ of (\ref
{eq1}) can be approximated by the solution $v$ of (\ref{eq1bis}) in the
following sense:
%
%
\begin{equation}
\label{finale}
\lim_{\varepsilon\to0}\E\sup_{t \in[\delta,T]}\biggl|\int_D
|u_\varepsilon(t,x)-v(t)|^2 \mu(dx)\biggr|^p=0
\end{equation}
for any fixed $0<\delta<T$ and $p\geq1/2$.

In order to prove (\ref{finale}), we first have to prove that for any
fixed $\varepsilon>0$, equation (\ref{eq1}) admits a unique adapted
mild solution in $L^p(\Omega,C([0,T];L^2(D)))$, that is, there exists
a unique adapted process $u_\varepsilon$ such that
\begin{eqnarray*}
u_\varepsilon(t)&=&e^{t A/\varepsilon}u_0+\int_0^te^{(t-s)
A/\varepsilon}F(s,u_\varepsilon(s)) \,ds+\int_0^te^{(t-s) A/\varepsilon
}G(s,u_\varepsilon(s)) \,dw^Q(s)\\
&&{}+w^\varepsilon_{A,B}(t),
\end{eqnarray*}
where $w^\varepsilon_{A,B}(t)$ is the boundary term (the \textit{stochastic
boundary convolution})
\[
w^\varepsilon_{A,B}(t)=(\delta_0-A)\int_0^te^{(t-s) A/\varepsilon
}N_{\delta
_0} [\Sigma(s)\,dw^B(s)],\qquad t\geq0
\]
(here, and in what follows, $F$ and $G$ denote the
composition/multiplication operators associated with $f$ and $g$, resp.).
In particular, we have to show that the above term is well defined in
$L^p(\Omega,C([0,T];L^2(D)))$. Concerning the notion of mild solutions
and existence and uniqueness results for SPDEs like (\ref{eq1}), with
fixed $\varepsilon>0$, we refer to Da Prato and Zabczyk
\cite{dapratozab}. However, we would like to stress that in the present
paper, we are not imposing the Hilbert--Schmidt condition on the
covariance operators $Q$ and $B$, and this makes the treatment of the
stochastic convolution and of the stochastic boundary convolution more
complicated, in view also of the a priori estimates with respect to
$\varepsilon>0$.

Actually, once we have a unique adapted mild solution $u_\varepsilon$
for (\ref{eq1}), we prove an a priori estimate of the
following type:
\[
\sup_{\varepsilon\in(0,1]}\E|u_\varepsilon
(t)|_{C([0,T];L^2(D))}^p\leq c_{T,p}\bigl(1+|u_0|_{L^2(D)}^p\bigr).
\]
Due to (\ref{dubbio}), this allows us to proceed to the proof of
(\ref{finale}).

After we have proven (\ref{finale}), in the final section, we study
the fluctuations of $u_\varepsilon$ from $v$. Namely, we introduce the
random field
\[
z_\varepsilon(t,x):=\frac{u_\varepsilon(t,x)-v(t)}{\sqrt{\varepsilon}},\qquad
(t,x) \in[0,+\infty)\times D,
\]
and show that, under the assumption that the noisy perturbation in
(\ref{eq1}) is of additive type (i.e., the diffusion
coefficient $g$ is independent of $u$), for any $t>0$,
\[
z_\varepsilon(t)\rightharpoonup I_0(t) \qquad\mbox{in } L^2(D,\mu),
\varepsilon\downarrow0,
\]
where $I_0(t,x)$ is the Gaussian random field taking values in
$L^2(D,\mu)$ for any $t>0$, defined by
\begin{eqnarray*}
I_0(t,x)&:=&\int_0^\infty\bigl(e^{sA}G(t)\,dw^Q(s,x)-\langle\hat
{G}(t),dw^Q(s)\rangle_{L^2(D)}\bigr)\\
&&{} + \int_0^\infty\bigl((\delta_0-A)e^{sA}N_{\delta_0}[\Sigma
(t)\,dw^B(s)](x)-\langle\hat{\Sigma}(t),
dw^B(s)\rangle_{L^2(\partial D)}\bigr).
\end{eqnarray*}

The random field $I_0(t,x)$ is well defined in $L^2(D,\mu)$ because of
the spectral gap inequality (\ref{dubbio}) and, in the case where the
coefficients $g$ and $\sigma$ do not depend on $t$, $I_0(t,x)$ also does
not depend on $t$ so that the weak limit of $z_\varepsilon(t,x)$ as
$\varepsilon\downarrow0$ depends only on the space variable $x$ and is
constant in time for any $t>0$.

\section{Notation and assumptions}
\label{sec2}

Let $D$ be a bounded domain in $\reals^d$, with \mbox{$d\geq1$}, satisfying
the extension and exterior cone properties, and let $\nu$ be the
outward normal at $\partial D$. We assume that $\partial D$ is a
$C^\infty$ manifold and $D$ is locally only on one side of $\partial
D$. In the case $d=1$, $D$ is a bounded interval $(a,b)$.

We define $H:=L^2(D)$ and $Z:=L^2(\partial D)$ and,
for any $\alpha\geq0$, we define $H^{\alpha}:=H^{\alpha}(D)$ and
$Z^{\alpha
}:=H^{\alpha}(\partial D)$ (in particular, $H^0=H$ and $Z^0=Z$).

We assume that $\mathcal{A}$ is a second order differential operator,
\[
\mathcal{A} f=\sum_{i,j=1}^d \frac{\partial}{\partial x_i}
\biggl(a_{ij}(x)\,\frac{\partial f}{\partial x_j}\biggr)+\sum_{i=1}^d
b_i(x)\,\frac{\partial f}{\partial x_i},
\]
satisfying the uniform ellipticity condition
\[
\inf_{x \in D} \sum_{i,j=1}^d a_{ij}(x)\xi_i \xi_j\geq a_0\sum
_{i=1}^d\xi_i^2,\qquad\xi\in\reals^d,
\]
for some $a_0>0$. The coefficients $a_{ij}$ and $b_i$ are assumed to be
smooth [for simplicity, we take them to be in $C^\infty(D)$]. In what
follows, we shall denote by $A$ the realization in $H$ of the operator
$\mathcal{A}$, endowed with the boundary condition
%
%
\begin{equation}
\label{boundary}
\frac{\partial h}{\partial\nu}(x):=\langle a(x)\nu(x), \nabla h(x)
\rangle_{\reals^{d}}=0,\qquad x \in\partial D.
\end{equation}
Namely,
\[
\cases{
A h=\mathcal{A} h, \qquad h \in D(A),\cr
D(A)=\{ h \in H^2(D) ; \langle a(x)\nu(x), \nabla h(x)
\rangle_{\reals^{d}}=0, x \in\partial D \}.}
\]
As is well known, the operator $A$ generates an analytic semigroup $\{
e^{tA}\}_{t\geq0}$ in $H$ which is also strongly continuous. Moreover,
\[
D(A^\alpha)=D((A^\star)^{\alpha})\subset H^{2\alpha}, \qquad\alpha\geq0,
\]
and
%
%
\begin{equation}
\label{dominio}
D(A^\alpha)=H^{2\alpha},\qquad 0\leq\alpha<\tfrac34
\end{equation}
(for proofs, see \cite{triebel} and \cite{lasiecka}, resp.).

If, for any $1<p\leq\infty$, we denote by $A_p$ the realization in
$L^p(D)$ of the operator~$\mathcal{A}$, endowed with the boundary condition
(\ref{boundary}), it can be proven that $A_p$
generates a strongly continuous analytic semigroup $e^{tA_p}$ in
$L^p(D)$. Notice that all of these semigroups are consistent, so, in
what follows, we shall denote them all by $e^{tA}$.

As proved in, for example, \cite{dav}, Theorem 2.4.4, since $\mathcal
{A}$ is uniformly elliptic and the domain $D$ has the extension
property, the semigroup $e^{tA}$ admits an integral kernel $k_t(x,y)$.
Due to the boundary condition, the kernel satisfies
%
%
\begin{equation}
\label{kernel}
0\leq k_t(x,y)\leq c (t^{-d/2}+ 1),\qquad t>0,
\end{equation}
for some constant $c>0$, almost everywhere in $D\times D$.

As a consequence of our assumptions on $\mathcal{A}$ and $D$, it is
possible to prove that there exists some $\delta_0 \in\reals$ such
that for any $\delta\geq\delta_0$ and $h \in Z$, the elliptic boundary
value problem
%
%
\begin{equation}
\label{neumann}
\cases{
(\delta-\mathcal{A})v(x)=0, &\quad $x \in D$,\cr
\langle a(x)\nu(x), \nabla v(x)\rangle_{\reals^{d}}=h(x),&\quad $x \in
\partial D$,}
\end{equation}
admits a unique weak solution $v \in H$, which we will denote by
$N_\delta h$. The application
$N_\delta\dvtx Z\to H$ is known as the \textit{Neumann map} associated
with the
operator $\mathcal{A}$.
It is well known that $N_\delta$ maps $Z$ into $H$ as a bounded linear
mapping. Moreover,
according to elliptic theory for domains with smooth boundaries (for a
proof, see~\cite{lionsmagenes}, Theorem 7.4 of Volume I), we have
%
%
\begin{equation}
\label{sobolev}
N_\delta\in\mathcal{L}(Z^{\alpha},H^{\alpha+3/2}),\qquad
\alpha\geq0.
\end{equation}

In what follows, we shall assume that $e^{tA}$ has the following
long-time behavior.
\begin{Hypothesis}
\label{H1}
The semigroup $e^{tA}$, $t\geq0$, admits a unique invariant measure
$\mu$ and there exists some $\gamma>0$ such that, for any $h \in
L^2(D,\mu)$,
%
%
\begin{equation}
\label{ergo}
\biggl|e^{tA} h-\int_D h(y)\mu(dy)\biggr|_{L^2(D,\mu)}\leq c
e^{-\gamma t} |h|_{L^2(D,\mu)},\qquad
t\geq0.
\end{equation}
\end{Hypothesis}

In what follows, we shall set $H_\mu:=L^2(D,\mu)$ and
\[
\langle h,\mu\rangle:=\int_D h(x)\mu(dx).
\]
\begin{Remark}
\label{rem21}
\begin{enumerate}
\item If $\mathcal{A}$ is a divergence-type operator, that is,
$b_i\equiv0$ for any $i=1,\ldots,d$, then the operator $A$ is
self-adjoint in $H$. This implies that it is possible to fix a complete
orthonormal system $\{e_k\}_{k\geq0}$ in $H$ and an increasing
sequence of nonnegative real numbers $\{\alpha_k\}_{k\geq0}$ such that
\[
A e_k=-\alpha_k e_k,\qquad
k \in\nat.
\]
Let $e_0$ be the constant eigenfunction corresponding to the $\alpha_0=0$
eigenvalue and let $\alpha_1$ be the first positive eigenvalue. It is
immediate to check that
%
%
\begin{equation}
\label{diag}
\mu(dx)=e_0^2 \,dx=|D|^{-1} \,dx
\end{equation}
and, in particular, that $H=H_\mu$, with equivalence of norms.
Moreover, as for any $x \in H$, we have
\[
e^{tA} x-\langle x,\mu\rangle=\sum_{i=1}^\infty e^{-t\alpha_i}\langle x,e_i\rangle_H e_i
\]
and $\alpha_1\leq\alpha_i$ for any $i\geq1$, it is immediate to check that
\[
|e^{tA} x-\langle x,\mu\rangle|^2_{H_\mu}=|D|^{-1} \sum_{i=1}^\infty
e^{-2 t\alpha_i}\langle x,e_i\rangle_H^2 \leq e^{-2 t\alpha_1} |x|_{H_\mu}^2,
\]
so the constant $\gamma$ in (\ref{ergo}) coincides with $\alpha_1$.

\item If $A$ is self-adjoint, as above, for any $\delta>0$ and $k \in
\nat$ it holds
that
%
%
\begin{equation}
\label{adjoint}
N^\star_\delta e_k=\frac1{\delta+\alpha_k} e_{k_{|_{\partial D}}}.
\end{equation}
Actually, for any $h \in Z$, we have
\begin{eqnarray*}
\langle N_\delta h,e_k\rangle_H&=&\frac1{\delta+\alpha_k}\int_D N_\delta h(x)(\delta
+\alpha
_k)e_k(x) \,dx\\
&=&\frac1{\delta+\alpha_k}\int_D N_\delta h(x)(\delta-\mathcal
{A})e_k(x) \,dx.
\end{eqnarray*}
Now, if we assume that $h \in Z^{1/2}$, according to (\ref
{sobolev}), we have that $N_\delta h \in H^2$ and then, due to the
Gauss--Green formula and to (\ref{neumann}), we obtain
\[
\int_D N_\delta h(x)\mathcal{A}e_k(x) \,dx=-\int_{\partial D} h(\sigma
) e_k(\sigma) \,d\sigma+\int_D \mathcal{A}N_\delta h(x)e_k(x) \,dx.
\]
This implies that
\begin{eqnarray*}
\langle N_\delta h,e_k\rangle_H &=& \frac1{\delta+\alpha_k}\int_D(\delta-\mathcal
{A})N_\delta h(x) e_k(x) \,dx+\frac1{\delta+\alpha_k}\int_{\partial D}
h(\sigma)
e_k(\sigma) \,d\sigma\\
&=&
\frac1{\delta+\alpha_k}\langle h,e_{k_{|_{\partial D}}}\rangle_Z
\end{eqnarray*}
so that
\[
\langle h,N^\star_\delta e_k\rangle_Z=\frac1{\delta+\alpha_k}\langle h,e_{k_{|_{\partial
D}}}\rangle_Z.
\]
As $Z^{1/2}$ is dense in $Z$, we can conclude that (\ref{adjoint}) holds.

\item As
\[
e^{tA}h(x)=\int_D k_t(x,y)h(y) \,dy,\qquad
x \in D,
\]
and $e^{tA}1=1$, we have
\[
|e^{tA}h(x)|^2\leq e^{tA}|h|^2(x),\qquad
x \in D.
\]
Due to the invariance of $\mu$, this implies that for any $h \in
H_\mu$,
\[
\int_D|e^{tA}h(x)|^2 \mu(dx)\leq\int_D e^{tA}|h|^2(x) \mu
(dx)=\int_D |h(x)|^2 \mu(dx),
\]
so $e^{tA}$ acts on $H_\mu$ as a contraction, that is,
%
%
\begin{equation}
\label{mar10}
\|e^{tA}\|_{{\mathcal L}(H_\mu)}\leq1,\qquad
t\geq0.
\end{equation}

\item We have that $H$ is continuously embedded into $H_\mu$.
Actually, due to the invariance of $\mu$ and to the kernel
representation of $e^{tA}$, for any $h \in H$, we have
\[
\int_D|h(x)|^2 \mu(dx)=\int_De^{1A}|h|^2(x) \mu(dx)=\int_D\int
_Dk_1(x,y)|h(y)|^2 \,dy \mu(dx).
\]
Then, thanks to (\ref{kernel}), we have
\[
|h|_{H_\mu}^2=\int_D|h(x)|^2 \mu(dx)\leq c \int_D|h(y)|^2 dy=|h|_H^2.
\]

\item As a matter of fact, there exists a nonnegative function $m \in
L^\infty(D)$ such that
\[
\mu(dx)=m(x) \,dx,\qquad
x \in D.
\]
Actually, let $\varphi, \psi\in C^2(\bar{D})$, with $\varphi$
fulfilling the boundary condition (\ref{boundary}).
Integrating by parts, we obtain
\[
\langle\psi,\mathcal{A}\varphi\rangle_{H}=\langle\mathcal{A}^\star\psi
,\varphi\rangle_H-\int_{\partial D} \langle a\nu, \nabla\psi\rangle_{\reals^{d}}
\varphi \,d\sigma+\int_{\partial D}\langle b,
\nu\rangle_{\reals^{d}}\varphi\psi\,
d\sigma,
\]
where
\[
\mathcal{A}^\star\psi=\frac{\partial}{\partial x_j}
\biggl(a_{ij}\,\frac{\partial\psi}{\partial x_i}\biggr)-\langle b, \nabla\psi
\rangle_{\reals^{d}}-\operatorname{div}b \psi.
\]
Hence, the operator $\mathcal{A}^\star$, endowed with the boundary condition
%
%
\begin{equation}
\label{adj}
\langle a(x)\nu(x),\nabla\psi(x)\rangle_{\reals^{d}} -\langle b(x), \nu(x)
\rangle_{\reals^{d}}\psi(x)=0,\qquad
x \in\partial D,
\end{equation}
is the formal adjoint of the operator $\mathcal{A}$, endowed with the
boundary condition (\ref{boundary}).

Now, the function $u=1$ is a nonzero solution of the problem
\[
\cases{
\mathcal{A}u(x)=0, &\quad
$x \in D$,\cr
\langle a(x)\nu(x), \nabla u(x)\rangle_{\reals^{d}}=0, &\quad
$x \in
\partial D$.}
\]
Then, by the Fredholm alternative, there exists a nonzero weak solution
$\varphi\in H^1$ to the adjoint problem
\[
\cases{
\mathcal{A}^\star\varphi(x)=0, &\quad
$x \in D$,\cr
\langle a(x)\nu(x), \nabla\varphi(x)\rangle_{\reals^{d}} -\langle b(x), \nu
(x)\rangle_{\reals^{d}}\varphi(x)=0, &\quad $x \in\partial D$.}
\]
By elliptic regularity results (cf. \cite{lu}, Chapter 3), as the
boundary of $D$ and the coefficients of $\mathcal{A}$ (and hence of
$\mathcal{A}^\star$) are of class $C^\infty$, we have that $\varphi
$ is a classical solution to the adjoint problem. Hence, if $A^\star$
is the adjoint of $A$ in $H$, for any $\lambda$ sufficiently large, we have
\[
(\lambda I-A^\star)^{-1}\varphi=\frac1\lambda\varphi
\]
and by taking the inverse Laplace transform, we obtain $e^{tA^\star
}\varphi=\varphi$ for any $t\geq0$.

Now, due to the positivity of the semigroup $e^{tA}$ (and hence of the
semigroup $e^{tA^\star}$) and to the fact that $e^{tA}$ is
conservative, we have that the set
\[
\Lambda:=\{ \varphi\in H \dvtx e^{tA^\star}\varphi=\varphi,
t\geq0 \}
\]
is a lattice, that is, $|\varphi| \in\Lambda$
for any $\varphi\in\Lambda$. Therefore, if we set
\[
m(x):=\frac{|\varphi(x)|}{\int_D|\varphi(y)| \,dy},\qquad x \in D,
\]
we have that $e^{tA^\star}m=m$ for any $t\geq0$ and hence
$m(x) \,dx$ is a probability measure and is invariant for $e^{tA}$. As
$\mu$ is the unique invariant measure for $e^{tA}$, we are done.
\end{enumerate}
\end{Remark}

Concerning the coefficients $f$, $g$ and $\sigma$ we assume the following
conditions.
\begin{Hypothesis}
\label{H2}
\begin{enumerate}
\item
The mappings $f, g\dvtx[0,\infty)\times D\times\reals\to\reals$ are
measurable and the mappings $f(t,x,\cdot), g(t,x,\cdot)\dvtx\reals\to
\reals$ are Lipschitz continuous, uniformly with respect to $(t, x)
\in[0,T]\times D$, for any $T>0$. Namely, for any $\xi,\eta\in
\reals$
\begin{eqnarray*}
{\sup_{(t, x) \in[0,T]\times D}}|f(t,x,\xi)-f(t,x,\eta)|&\leq&
L_{T,f} |\xi-\eta|,\\
{\sup_{(t, x) \in[0,T]\times D}}|g(t,x,\xi)-g(t,x,\eta)|&\leq&
L_{T,g} |\xi-\eta|.
\end{eqnarray*}

\item The mapping $\sigma\dvtx[0,\infty)\times\partial D \to\reals$ is
measurable
and for any $T>0$,
\[
{\sup_{t \in[0,T]}}|\sigma(t,\cdot)|_{L^\infty(\partial D)}=:c_{T,\sigma
}<\infty.
\]
\end{enumerate}
\end{Hypothesis}

In what follows, for any $t\geq0$ and $h_1,h_2 \in H$, we shall define
\[
F(t,h_1)(x):=f(t,x,h_1(x)),\qquad x \in D,
\]
and
\[
[G(t,h_1)h_2](x):=g(t,x,h_1(x))h_2(x),\qquad x \in D.
\]
Due to Hypothesis \ref{H2}, we have that
$F(t,\cdot)\dvtx H\to H$, $G(t,\cdot)\dvtx H\to{\mathcal L}(H,L^1(D))$
and $G(t,\cdot
)\dvtx H\to{\mathcal L}(L^\infty(D),H)$ are all Lipschitz continuous, uniformly
with respect to $t \in[0,T]$, for any $T>0$.

Notice that the same is true for the mappings
$F(t,\cdot)\dvtx H_\mu\to H_\mu$, $G(t,\cdot)\dvtx H_\mu\to{\mathcal
L}(H_\mu
,L^1(D,\mu))$ and $G(t,\cdot)\dvtx H_\mu\to{\mathcal L}(L^\infty(D;\mu
),H_\mu)$.

Analogously, if, for any $t\geq0$ and $z \in Z$, we set
\[
[\Sigma(t)z](x):=\sigma(t,x)z(x),\qquad x \in\partial D,
\]
then we have that $\Sigma(t)$ is a bounded linear operator on $Z$ and
for any $T>0$,
%
%
\begin{equation}
\label{sigma}
\|\Sigma(t)\|_{{\mathcal L}(Z)}\leq c_{T,\sigma},\qquad t \in[0,T].
\end{equation}

Finally, concerning the noisy perturbations $w^Q(t)$ and $w^B(t)$, we
assume that they are two independent cylindrical Wiener processes
defined on the same stochastic basis $(\Omega,\mathcal{F},\mathcal
{F}_t,\Pro)$, taking values in $H$ and $Z$, respectively, with
respective covariance operators $Q \in\mathcal{L}^+(H)$ and $B \in
\mathcal{L}^+(Z)$. Namely,
\[
w^Q(t)=\sum_{k \in\nat} \lambda_k e_k \beta_k(t),\qquad
w^B(t)=\sum_{k \in\nat} \theta_k f_k \hat{\beta}_k(t),
\]
where $\{e_k\}_{k \in\nat}$ is the orthonormal basis of $H$ which
diagonalizes $Q$, with eigenvalues $\{\lambda_k\}_{k \in\nat}$, $\{
f_k\}_{k \in\nat}$ is the orthonormal basis of $Z$ which
diagonalizes $B$, with eigenvalues $\{\theta_k\}_{k \in\nat}$, and
$\{\beta_k\}_{k \in\nat}$ and $\{\hat{\beta}_k\}_{k \in\nat
}$ are two sequences of independent standard Brownian motions, both
defined on the stochastic basis $(\Omega,\mathcal{F},\mathcal
{F}_t,\Pro)$. Notice that the two sequences above are not convergent
in $H$ and $Z$, but in any Hilbert spaces $U$ and $V$ which contain $H$
and $Z$, respectively, with Hilbert--Schmidt embedding. Moreover, in
the case $d=1$, we have $Z=\reals^2$ and hence
\[
w^B(t)=\Theta \hat{\beta}(t),
\]
where $\Theta=\operatorname{diag}(\theta_1,\theta_2) $ and $\hat
{\beta}(t)=(\hat{\beta}_1(t),\hat{\beta}_2(t))$ is a
two-dimensional standard Brownian motion.

In what follows, we shall assume the following summability conditions
on the eigenvalues $\lambda_k$ and $\theta_k$ and the sup-norm of the
corresponding eigenfunctions.
\begin{Hypothesis}
\label{H3}
\begin{enumerate}
\item If $d\geq2$, then there exists $\rho<2d/(d-2)$ such that
%
%
\begin{equation}
\label{lambda}
\sum_{k \in\nat} \lambda_k^\rho|e_k|_\infty^2=:\kappa_Q<\infty.
\end{equation}

\item If $d\geq2$, then there exists $\beta<2d/(d-1)$ such that
%
%
\begin{equation}
\label{theta}
\sum_{k \in\nat} \theta_k^\beta=:\kappa_B<\infty.
\end{equation}
\end{enumerate}
\end{Hypothesis}
\begin{Remark}
\label{Rem22}
\begin{enumerate}

\item From the proofs of Lemmas \ref{lemma33}, \ref{mar4} and
\ref{tax10}, it is possible to see that if the mapping
$g\dvtx[0,T]\times D\times\reals\to\reals$ is uniformly bounded for any
$T>0$, then we do not need to require that the sequence $\{e_k\}_{k \in
\nat}$ is contained in $L^\infty(D)$ and condition (\ref{lambda})
can be replaced by
\[
\sum_{k \in\nat} \lambda_k^\rho<\infty.
\]

\item As both $d/(d-2)$ and $d/(d-1)$ are strictly greater than $1$,
neither $Q$ nor $B$ are required to be Hilbert--Schmidt operators in
general. Moreover, in space dimension $d=1$, we have no conditions on
the eigenvalues $\{\lambda_k\}$ and
we can take $Q=I$. This means that we can deal with space--time white noise.
\end{enumerate}
\end{Remark}

\section[A priori bounds for the solution of (1.1)]{A priori bounds for the solution of (\protect\ref{eq1})}
\label{sec3}

In this section, we are concerned with uniform bounds for the $p$th
moments of the $C([0,T];H)$-norm of the mild solution $u_\varepsilon$ of
(\ref{eq1}).

We first recall some general facts about the linear parabolic equation
with nonhomogeneous boundary conditions
%
%
\begin{equation}
\label{detlin}
\cases{
\displaystyle\frac{\partial y}{\partial
t}(t,x)=\mathcal{A}y(t,x),&\quad
$t\geq0, x \in D$,\vspace*{2pt}\cr
\langle a(x)\nu(x),\nabla y(t,x)\rangle_{\reals^{d}}=v(t,x),&\quad $t\geq0,
x \in\partial D$,\vspace*{2pt}\cr
y(0,x)=y_0(x), &\quad $x \in D$,}
\end{equation}
where $v$ is a $Z$-valued function. If $v(\cdot)$ is twice
continuously differentiable and there exists $\delta_0>0$ such that
$y_0-N_\delta v(0) \in D(A)$ for $\delta>\delta_0$, then the solution of
problem (\ref{detlin})
is given by
%
%
\begin{equation}
\label{solution}
y(t)=e^{tA}y_0+(\delta-A)\int_0^te^{(t-s)A}N_\delta v(s) \,ds
\end{equation}
(for a proof, see, e.g., \cite{dpz2}, Proposition 13.2.1).

Such a formula can be extended by continuity to less regular functions $v$.
In particular, for each $\varepsilon>0$, we can consider the problem
%
%
\begin{equation}
\label{stoclin}\quad
\cases{
\displaystyle\frac{\partial y}{\partial t}(t,x)=\frac1\varepsilon\mathcal
{A}y(t,x), &\quad $t\geq0, x \in D$,\vspace*{2pt}\cr
\displaystyle\langle a(x)\nu(x), \nabla y(t,x)\rangle_{\reals^{d}}=\varepsilon\sigma
(t,x)\,\frac{\partial w^B}{\partial t}(t,x), &\quad $t\geq0, x \in
\partial D$,\vspace*{2pt}\cr
y(0,x)=0, &\quad $x \in D$,}
\end{equation}
where $w^B$ is the cylindrical Wiener process defined in $Z$,
introduced in Section \ref{sec2}. In analogy to formula (\ref
{solution}), by taking $\delta=\delta_0/\varepsilon$ and
$v(t)=\varepsilon
\Sigma(t) \,\partial w^B/\partial t$, we say that for any $\varepsilon\in
(0,1]$, the process
\[
w^\varepsilon_{A,B}(t)=(\delta_0-A)\int_0^te^{(t-s) A/\varepsilon
}N_{\delta
_0} [\Sigma(s)\,dw^B(s)],\qquad t\geq0,
\]
is a \textit{mild solution} to problem (\ref{stoclin}).
The process $w^\varepsilon_{A,B}(t)$ can be interpreted as a
\textit{boundary Ornstein--Uhlenbeck} process and can be written as the
infinite series
\[
w^\varepsilon_{A,B}(t)=\sum_{k \in\nat} (\delta_0-A)\int
_0^te^{(t-s) A/\varepsilon}N_{\delta_0} [\Sigma(s)Bf_k] \,d\hat
{\beta}_k(s),\qquad t\geq0.
\]
As proved in the next lemma, such a series is well defined in
$L^p(\Omega;C([0,T];H))$ for any $T>0$ and $p\geq1$. Moreover, a
uniform estimate with respect to $\varepsilon\in(0,1]$ holds.
\begin{Lemma}
\label{lemma31}
Under part 2 of Hypothesis \ref{H3}, the process $w^\varepsilon_{A,B}$
belongs to $L^p(\Omega;C([0,T];H))$ for any $T>0$, $p\geq1$ and
$\varepsilon\in(0,1]$,
and
%
%
\begin{equation}
\label{aprioriB}
\sup_{\varepsilon\in(0,1]}\E|w^\varepsilon
_{A,B}|^p_{C([0,T];H)}=:c_{T,p}<\infty.
\end{equation}
\end{Lemma}
\begin{pf}
As a consequence of the stochastic Fubini theorem and of the elementary identity
\[
\int_\sigma^t(t-s)^{\alpha-1}(s-\sigma)^{-\alpha} \,ds=\frac{\pi}{\sin\pi
\alpha},\qquad 0\leq\sigma\leq t, \alpha\in(0,1),
\]
we have the
factorization formula
\[
w^\varepsilon_{A,B}(t)=\frac{\sin\pi\alpha}{\pi}\int_0^t(t-s)^{\alpha
-1}e^{(t-s)A/\varepsilon}Y_{\varepsilon,\alpha}(s) \,ds,
\]
where
\[
Y_{\varepsilon,\alpha}(s)=\int_0^s(s-r)^{-\alpha}(\delta
_0-A)e^{(s-r)A/\varepsilon}N_{\delta_0} [\Sigma(r)\,dw^B(r)]
\]
(for a proof, see \cite{dpkz}).
By the H\"older inequality, this implies that for any $\alpha>1/p$,
%
%
\begin{eqnarray}
\label{1}\hspace*{18pt}
&&
{\E\sup_{t \in[0,T]}} |w^\varepsilon_{A,B}(t)|_H^p\nonumber\\
&&\qquad\leq
c_{T,p,\alpha}\int_0^T \E|Y_{\varepsilon,\alpha}(s)|_H^p
\,ds\nonumber\\[-8pt]\\[-8pt]
&&\qquad\leq c_{T,p,\alpha}\int_0^T\biggl(\int_0^s(s-r)^{-2\alpha}\nonumber\\
&&\qquad\quad\hspace*{63.4pt}{}\times\sum_{k \in
\nat}\theta_k^2 \bigl|(\delta_0-A)e^{(s-r)A/\varepsilon}N_{\delta_0}
[\Sigma(r)f_k]\bigr|_H^2 \,dr\biggr)^{p/2} \,ds,\nonumber
\end{eqnarray}
the last inequality following from the Burkholder--Davis--Gundy inequality.

Now, assume that $d>1$ (the case $d=1$ is simpler). According to (\ref
{theta}), we have
%
%
\begin{eqnarray}
\label{2}
&&
\sum_{k \in\nat}\theta_k^2 \bigl|(\delta_0-A)e^{(s-r)
A/\varepsilon}N_{\delta_0} [\Sigma(r)f_k]\bigr|_H^2\nonumber\\
&&\qquad\leq\kappa_B^{2/\beta}\biggl(\sum_{k \in\nat}\bigl|(\delta
_0-A)e^{(s-r)A/\varepsilon}N_{\delta_0} [\Sigma(r)f_k]
\bigr|_H^2\biggr)^{1/\zeta}\\
&&\qquad\quad\hspace*{0pt}{}\times
\sup_{k \in\nat}\bigl|(\delta_0-A)e^{(s-r)A/\varepsilon}N_{\delta_0}
[\Sigma(r)f_k]\bigr|_H^{{2(\zeta-1)}/\zeta},\nonumber
\end{eqnarray}
where
$\zeta:=\beta/(\beta-2)$. Thanks to (\ref{dominio}) and (\ref
{sobolev}), for any $\rho>0$, we have
%
%
\begin{equation}
\label{schauder}
S_\rho:=(\delta_0-A)^{({3-\rho})/4}N_{\delta_0} \in\mathcal{L}(Z,H).
\end{equation}
Hence, for any $\varepsilon>0$ and $0\leq r\leq s\leq T$, due to (\ref
{sigma}), we have
%
%
\begin{eqnarray}
\label{somma}\hspace*{22pt}
&&
\sum_{k \in\nat}\bigl|(\delta_0-A)e^{(s-r)A/\varepsilon}N_{\delta
_0} [\Sigma(r)f_k]\bigr|_H^2\nonumber\\
&&\qquad=
\sum_{k \in\nat}\bigl|e^{{(s-r)}/2 A/\varepsilon}(\delta
_0-A)^{({1+\rho})/4}e^{{(s-r)}/2 A/\varepsilon}S_\rho
\Sigma(r)f_k\bigr|_H^2\nonumber\\
&&\qquad=\sum_{k \in\nat} \sum_{h \in\nat}\bigl|\bigl\langle f_k,\Sigma
(r)S_\rho^\star \bigl[(\delta_0-A)^{({1+\rho})/ 4}e^{
{(s-r)}/2 A/\varepsilon}\bigr]^\star e^{{(s-r)}/2{A^\star
}/\varepsilon} e_h\bigr\rangle_Z\bigr|^2\\
&&\qquad=\sum_{h \in\nat}\bigl|\Sigma(r)S_\rho^\star \bigl[(\delta
_0-A)^{({1+\rho})/4}e^{{(s-r)}/2 A/\varepsilon}\bigr]^\star
e^{{(s-r)}/2{A^\star}/\varepsilon} e_h\bigr|_Z^2\nonumber\\
&&\qquad\leq c_{T,\rho} \biggl[\biggl(\frac\varepsilon{s-r}\biggr)^{
({1+\rho})/2}+1\biggr]\sum_{h \in\nat}\bigl|e^{{(s-r)}/2
{A^\star}/\varepsilon} e_h\bigr|_H^2.\nonumber
\end{eqnarray}

As the semigroup $e^{tA}$ admits an integral kernel $k_t(x,y)$,
that is,
\[
e^{tA}f(x)=\int_D k_t(x,y)f(y) \,dy,\qquad x \in D,
\]
we have
\[
e^{tA^\star} h(y)=\int_Dk_t(x,y)h(x) \,dx,\qquad y \in D.
\]
This implies
%
%
\begin{eqnarray}
\label{sommabis}
\sum_{h \in\nat}\bigl|e^{{(s-r)}/2{A^\star
}/\varepsilon} e_h\bigr|_H^2&=&
\sum_{h \in\nat}\int_D\bigl|e^{{(s-r)}/2{A^\star
}/\varepsilon} e_h(y)\bigr|^2 \,dy\nonumber\\
&=&\sum_{h \in\nat}\int_D\biggl|\int_D k_{
({s-r})/({2\varepsilon})}(x,y)e_h(x) \,dx\biggr|^2 \,dy\nonumber\\[-8pt]\\[-8pt]
&=&\sum_{h \in\nat}\int
_D\bigl|\bigl\langle k_{({s-r})/({2\varepsilon})}(\cdot,y),e_h\bigr\rangle_H\bigr|^2
\,dy\nonumber\\
&=&\int_D\bigl|k_{({s-r})/({2\varepsilon})}(\cdot,y)\bigr|_H^2 \,dy.\nonumber
\end{eqnarray}
Now, due to (\ref{kernel}), for any $t>0$ and $y \in D$, we have
\[
|k_{t}(\cdot,y)|_H^2=\int_D|k_t(x,y)|^2 \,dx\leq c(t^{-
d/2}+ 1)\int_Dk_t(x,y) \,dx
\]
and hence
\[
\int_D|k_{t}(\cdot,y)|_H^2 \,dy\leq c(t^{-d/2}+ 1
)\int_{D\times D}k_t(x,y) \,dx \,dy=c |D| (t^{-d/2}+ 1).
\]
This implies that for any $\varepsilon>0$,
\[
\sum_{h \in\nat}\bigl|e^{{(s-r)}/2{A^\star}/\varepsilon}
e_h\bigr|_H^2\leq
c |D| \biggl[\biggl(\frac\varepsilon{s-r}\biggr)^{d/2}+ 1\biggr],
\]
so, thanks to (\ref{somma}), we have
%
%
\begin{eqnarray}
\label{quadrato}
&&\biggl( \sum_{k \in\nat}\bigl|(\delta_0-A)e^{(s-r)A/\varepsilon
}N_{\delta_0}[\Sigma(r) f_k]\bigr|_H^2\biggr)^{1/\zeta}\nonumber\\[-8pt]\\[-8pt]
&&\qquad\leq
c_{T,\rho} \biggl[\biggl(\frac\varepsilon{s-r}\biggr)^{({d+1+\rho
})/({2\zeta})}+ 1\biggr].\nonumber
\end{eqnarray}

Next, by proceeding as in (\ref{somma}), we have
%
%
\begin{eqnarray}
\label{3}
&&\sup_{k \in\nat}\bigl|(\delta_0-A)e^{(s-r)A/\varepsilon}N_{\delta_0}
[\Sigma(r)f_k]\bigr|_H^{{2(\zeta-1)}/\zeta}\nonumber\\[-8pt]\\[-8pt]
&&\qquad\leq c_{T,\rho
} \biggl[\biggl(\frac\varepsilon{s-r}\biggr)^{{(1+\rho)(\zeta
-1)}/({2\zeta})}+1\biggr].\nonumber
\end{eqnarray}
Therefore, thanks to (\ref{1}), (\ref{2}), (\ref{quadrato}) and
(\ref{3}), we can conclude that for any $\varepsilon\in(0,1]$,
\[
\E\sup_{t \in[0,T]} |w^\varepsilon_{A,B}(t)|_H^p\leq c_{T,p,\alpha
,\rho}\biggl(\int_0^T\bigl[s^{-(2\alpha+({d+\zeta})/({2\zeta})+
{\rho}/2)}+ 1\bigr] \,ds\biggr)^{p/2}.
\]

Now, as in Hypothesis \ref{H3}, we are assuming that $\beta
<2d/(d-1)$, so we have $(d+\zeta)/2\zeta<1$. This means that we can fix
$\bar{\alpha}>0$ and $\bar{\rho}>0$ such that
\[
2\bar{\alpha}+\frac{d+\zeta}{2\zeta}+\frac{\bar{\rho}}2<1
\]
and then, for any $p>\bar{p}:=1/\bar{\alpha}$ we obtain
\[
\sup_{\varepsilon\in(0,1]} \E\sup_{t \in[0,T]} |w^\varepsilon
_{A,B}(t)|_H^p\leq c_{T,p}.
\]
The estimate for general $p\geq1$ follows from the H\"older inequality.
\end{pf}

Next, we pass to (\ref{eq1}).
\begin{Definition}
Let $T>0$ and $p\geq1$. An adapted process $u_\varepsilon\in
L^p(\Omega;C([0$,\break $T];H))$ is a \textit{mild} solution of (\ref
{eq1}) if, for any $t \in[0,T]$,
\[
u_\varepsilon(t)=e^{tA/\varepsilon}u_0+\int_0^te^{(t-s)
A/\varepsilon}F(s,u_\varepsilon(s)) \,ds+w^\varepsilon_{A,Q}(u_\varepsilon
)(t)+w^\varepsilon_{A,B}(t),
\]
where, for any $u \in L^p(\Omega;C([0,T];H))$, we define
\[
w^\varepsilon_{A,Q}(u)(t):=\int_0^te^{(t-s) A/\varepsilon}G(s,u(s))\,
dw^Q(s),\qquad t\geq0.
\]
As is well known, $w^\varepsilon_{A,Q}(u)$ is the unique mild solution
of the problem
%
%
\begin{equation}
\label{conv}\quad
\cases{
\displaystyle\frac{\partial y}{\partial t}(t,x)=\frac1\varepsilon\mathcal
{A}y(t,x)+g(t,x,u(t,x)) \,\frac{\partial w^Q}{\partial t}(t,x),&\quad $t\geq0, x
\in D$,\vspace*{2pt}\cr
\langle a(x)\nu(x), \nabla y(t,x)\rangle_{\reals^{d}}=0, &\quad $t\geq0,
x \in\partial D$,\vspace*{2pt}\cr
y(0,x)=0,&\quad $x \in D$,}\hspace*{-27pt}
\end{equation}
where $w^Q$ is the cylindrical Wiener process with values in $H$,
introduced in Section \ref{sec2}.

As for $w^\varepsilon_{A,B}$, we show that $w^\varepsilon_{A,Q}$
satisfies a bound in $L^p(\Omega;C([0,T];H))$ which is uniform with
respect to $\varepsilon\in(0,1]$.
\end{Definition}
\begin{Lemma}
\label{lemma33}
Assume Hypothesis \ref{H2} and part 1 of Hypothesis \ref{H3}. Then,
$w^\varepsilon_{A,Q}$ is Lipschitz continuous from $L^p(\Omega
;C([0,T];H))$ into itself for any $T>0$ and $p\geq1$,
and
%
%
\begin{equation}
\label{aprioriQ}
\sup_{\varepsilon\in(0,1]}\E|w^\varepsilon
_{A,Q}(u)|^p_{C([0,T];H)}\leq c_{T,p}\biggl(1+\E\int_0^T|u(s)|_H^p
\,ds\biggr).
\end{equation}
\end{Lemma}
\begin{pf}
The proof of the Lipschitz continuity of $w^\varepsilon_{A,Q}$ in
$L^p(\Omega;C([0,T]$; $H))$ is classical and can be found in, for
example, \cite{cerrai1}. Concerning estimate (\ref{aprioriQ}), as in
the proof of Lemma \ref{lemma31}, we use a factorization argument and,
for any $\alpha>1/p$, we get
\begin{eqnarray*}
&&\E\sup_{t \in[0,T]} |w^\varepsilon_{A,Q}(t)|_H^p\\
&&\qquad\leq c_{T,p,\alpha} \E\int_0^T\biggl(\int_0^s(s-r)^{-2\alpha}\sum_{k
\in\nat}\lambda_k^2 \bigl|e^{(s-r) A/\varepsilon}
[G(r,u(r))e_k]\bigr|_H^2 \,dr\biggr)^{p/2} \,ds.
\end{eqnarray*}
According to (\ref{lambda}), if we set $\zeta:=\rho/(\rho-2)$, then
we have
%
%
\begin{eqnarray}
\label{22}
&&\sum_{k \in\nat}\lambda_k^2 \bigl|e^{(s-r) A/\varepsilon}
[G(r,u(r))e_k]\bigr|_H^2\nonumber\\
&&\qquad\leq\kappa_Q^{2/\rho}\biggl(\sum_{k \in\nat}
\bigl|e^{(s-r) A/\varepsilon}[G(r,u(r))e_k]\bigr|_H^2\biggr)^{1/
\zeta}\\
&&\qquad\quad{}\times\sup_{k \in\nat}\bigl|e^{(s-r)A/\varepsilon}[G(r,u(r))e_k
]\bigr|_H^{{2(\zeta-1)}/\zeta}|e_k|_\infty^{-4/\rho}.\nonumber
\end{eqnarray}
As in the proof of (\ref{sommabis}), we have
\[
\sum_{k \in\nat}\bigl|e^{(s-r) A/\varepsilon}[G(r,u(r))e_k
]\bigr|_H^2=\int_D
\bigl|k_{({s-r})/\varepsilon}(x,\cdot)g(r,\cdot,u(r))\bigr|_H^2 \,dx.
\]
Now, thanks to (\ref{kernel}), for any $t>0$, $x \in D$ and $h \in H$,
we have
%
%
\begin{eqnarray}
\label{mar1}
|k_t(x,\cdot)h|_H^2&=&\int_D |k_t(x,y) h(y)|^2 \,dy\nonumber\\
&\leq& c (t^{-d/2}+1 )\int_D k_t(x,y)h^2(y) \,dy\\
&=&c (t^{-
d/2}+ 1)e^{tA}h^2(x)\nonumber
\end{eqnarray}
and this is meaningful since $e^{tA}$ is well defined in $L^1(D)$.
In particular, for any $\varepsilon>0$,
\begin{eqnarray*}
&&\sum_{k \in\nat}\bigl|e^{(s-r)A/\varepsilon}
[G(r,u(r))e_k]\bigr|_H^2\\
&&\qquad\leq c\biggl[\biggl(\frac\varepsilon{s-r}
\biggr)^{d/2}+ 1\biggr]\int_D e^{(s-r)A/\varepsilon}g^2(r,\cdot
,u(r))(x) \,dx\\
&&\qquad=c \biggl[\biggl(\frac\varepsilon{s-r}\biggr)^{d/2}+ 1
\biggr]\bigl|e^{(s-r)A/\varepsilon}g^2(r,\cdot,u(r))\bigr|_{L^1(D)}\\
&&\qquad\leq c
\biggl[\biggl(\frac\varepsilon{s-r}\biggr)^{d/2}+ 1\biggr]|g(r,\cdot,u(r))|_H^2
\end{eqnarray*}
and, due to the linear growth of $g$,
%
%
\begin{eqnarray}
\label{33}
&&\biggl(\sum_{k \in\nat}\bigl|e^{(s-r) A/\varepsilon}
[G(r,u(r))e_k]\bigr|_H^2\biggr)^{1/\zeta}\nonumber\\[-8pt]\\[-8pt]
&&\qquad\leq c_T \biggl[\biggl(\frac\varepsilon{s-r}\biggr)^{d/({2\zeta})}+ 1
\biggr]\bigl(1+|u(r)|_H^{2/\zeta}\bigr).\nonumber
\end{eqnarray}
By analogous arguments, we have
%
%
\begin{equation}
\label{33bis}\hspace*{28pt}
\sup_{k \in\nat}
\bigl|e^{(s-r) A/\varepsilon}[G(r,u(r))e_k]\bigr|_H^{
{2(\zeta-1)}/\zeta}|e_k|_\infty^{-4/\rho}\leq c_T
\bigl(1+|u(r)|_H^{{2(\zeta-1)}/\zeta}\bigr)
\end{equation}
and then, thanks to (\ref{22}), (\ref{33}) and (\ref{33bis}), we
get, for any $\varepsilon\in(0,1]$,
\begin{eqnarray*}
&&\E\sup_{t \in[0,T]} |w^\varepsilon_{A,Q}(t)|_H^p\\
&&\qquad\leq c_{T,p,\alpha} \E\int_0^T\biggl(\int_0^s\biggl[\biggl(\frac1
{s-r}\biggr)^{2\alpha+{d}/({2\zeta})}+1\biggr]\bigl(1+|u(r)|_H^2\bigr)\, dr
\biggr)^{p/2}\, ds.\nonumber
\end{eqnarray*}
As we are assuming $\rho<2d/(d-2)$, we can find $\bar{\alpha}>0$ such
that $2\bar{\alpha}+d/(2\zeta)<1$. Due to the Young inequality, this
implies (\ref{aprioriQ}) for all $p>\bar{p}=1/\bar{\alpha}$ and hence
for all $p\geq1$.
\end{pf}

According to Lemmas \ref{lemma31} and \ref{lemma33}, we have the
following result.
\begin{Theorem}
\label{sol}
Under Hypotheses \ref{H1}, \ref{H2} and \ref{H3}, for any $T>0$ and
$p\geq1$, and for any $u_0 \in H$ and $\varepsilon>0$, equation
(\ref{eq1}) admits a unique adapted mild solution $u_\varepsilon\in
L^p(\Omega;C([0,T];H))$. Moreover,
%
%
\begin{equation}
\label{epsilon}
\sup_{\varepsilon\in(0,1]}\E|u_\varepsilon|_{C([0,T];H)}^p\leq
c_{T,p}(1+|u_0|_H^p).
\end{equation}
\end{Theorem}
\begin{pf}
As both $F(t,\cdot)\dvtx H\to H$ and $w^\varepsilon_{A,Q}\dvtx
L^p(\Omega; C([0,T];H))\to L^p(\Omega$; $C([0,T];H))$ are Lipschitz continuous and
$w^\varepsilon_{A,B} \in L^p(\Omega;C([0,T];H))$, we have that the
mapping $\Phi_\varepsilon$ defined by
\[
\Phi_\varepsilon(u)(t)=e^{t A/\varepsilon}u_0+\int_0^t
e^{(t-s) A/\varepsilon}F(s,u(s)) \,ds+w^\varepsilon
_{A,Q}(u)(t)+w^\varepsilon_{A,B}(t)
\]
is Lipschitz continuous from the space of adapted processes in
$L^p(\Omega;C([0,T]$; $H))$ into itself. Therefore, by a classical fixed
point argument,
equation (\ref{eq1}) admits a unique adapted mild solution
$u_\varepsilon\in L^p(\Omega,C([0,T];H))$.

Next, for any $\varepsilon>0$, we have
\begin{eqnarray*}
|u_\varepsilon(t)|_H^p&\leq&c_p \biggl(|u_0|_H^p+c t^{p-1}\int_0^t
\bigl(1+|u_\varepsilon(s)|_H^p\bigr) \,ds\\
&&\hspace*{27.2pt}{}+|w^\varepsilon_{A,Q}(u_\varepsilon
)(t)|_H^p+|w^\varepsilon_{A,B}(t)|_H^p\biggr)
\end{eqnarray*}
and then, according to (\ref{aprioriB}) and (\ref{aprioriQ}), we
conclude that
\[
\E\sup_{t \in[0,T]} |u_\varepsilon(t)|_H^p\leq c_{T,p}
\biggl(1+|u_0|_H^p+\int_0^T\E\sup_{r \in[0,s]}|u_\varepsilon(r)|_H^p
\,ds\biggr).
\]
The Gronwall lemma allows us to obtain (\ref{epsilon}).
\end{pf}

\section{The averaging result}
\label{sec4}

In this section, we show that for any $0<\delta<T$ and $p\geq1$, the
sequence $\{u_\varepsilon\}_{\varepsilon\in(0,1]}$ converges in
$L^p(\Omega;C([\delta,T];H_\mu))$ to the solution of a suitable
one-dimensional stochastic differential equation. In what follows, we
first introduce the limiting equation by constructing the coefficients
and by describing a situation in which they are given by a simple
expression. In the second part of this section, we prove the
convergence result.

We start with the drift term.
For each $t\geq0$ and $h \in H$, we define
%
%
\begin{equation}
\label{fhat}
\hat{F}(t,h):=\langle F(t,h),\mu\rangle=\int_Df(t,x,h(x)) \mu(dx),
\end{equation}
where $\mu(dx)$ is the unique invariant measure associated with the
semigroup $e^{tA}$ (see Section \ref{sec2} and Hypothesis \ref{H1}).
According to Hypothesis \ref{H2}, for any $T>0$ and $h_1, h_2 \in
H$, we have
\[
|f(t,x,h_1(x))-f(t,x,h_2(x))|\leq L_{T,f} |h_1(x)-h_2(x)|,\qquad
(t,x) \in[0,T]\times D,
\]
so that
\[
\hat{F}(t,\cdot)\dvtx H_\mu\to\reals
\]
is Lipschitz continuous, uniformly with respect to $t \in[0,T]$, for
any $T>0$. Notice that, as $H\subset H_\mu$, this implies that $\hat
{F}(t,\cdot)\dvtx H\to\reals$ is also Lipschitz continuous.

Next, we construct the term arising from the stochastic
convolution\break
$w^\varepsilon_{A,Q}(u)(t)$. For each $t\geq0$ and $h \in H$, we
introduce the linear mapping
\[
z \in H\mapsto\sum_{k \in\nat}\langle G(t,h)e_k,\mu\rangle
\langle z,e_k\rangle_H=\langle G(t,h)z,\mu\rangle \in\reals.
\]
As $H$ is continuously embedded into $H_\mu$, for any $T>0$, we have
\[
|\langle G(t,h)z,\mu\rangle|\leq|g(t,\cdot,h)|_{H_\mu}|z|_{H_\mu
}\leq c_T(1+|h|_{H_{\mu}})|z|_H,\qquad t\leq T.
\]
This means that there exists $\hat{G}(t,h) \in H$ such that
\[
\langle\hat{G}(t,h),z\rangle_H=\langle G(t,h)z,\mu\rangle,\qquad z \in H.
\]
Moreover, since for any $h_1, h_2 \in H_\mu$ and $T>0$,
\begin{eqnarray*}
&&|\langle G(t,h_1)z,\mu\rangle-\langle G(t,h_2)z,\mu\rangle|\\
&&\qquad\leq
|g(t,\cdot,h_1)-g(t,\cdot,h_2)|_{H_\mu}|z|_{H}\\
&&\qquad\leq c_T |h_1-h_2|_{H_\mu}|z|_{H},\qquad t\leq T,
\end{eqnarray*}
we have that the mapping $\hat{G}(t,\cdot)\dvtx H_\mu\to H$ is Lipschitz
continuous, uniformly with respect to $t \in[0,T]$, for any $T>0$.

This, in particular, implies that the mapping $\hat{G}(t,\cdot)$ is
also Lipschitz continuous, both in $H$ and in $H_\mu$, uniformly for
$t \in[0,T]$.

Finally, we construct the term arising from the boundary convolution
$w^\varepsilon_{A,B}(t)$. For each fixed $t\geq0$, we introduce the mapping
\[
h \in Z\mapsto\delta_0\langle N_{\delta_0}[\Sigma(t)h],\mu\rangle=\delta_0\int_D
N_{\delta_0}[\sigma(t,\cdot)h](x) \mu(dx) \in\reals.
\]
As $N_{\delta_0}$ is a bounded linear operator from $Z$ into $H$,
$\Sigma
(t)$ is bounded and linear in $Z$ and $H$ is continuously embedded in
$H_\mu$, such a mapping is bounded and linear from $Z$ into $\reals$
and then,
for any $t\geq0$, there exists $\hat{\Sigma}(t) \in Z$ such that
for any $h \in Z$, we have
%
%
\begin{equation}
\label{sigmahat}
\langle\hat{\Sigma}(t),h\rangle_Z=\delta_0\langle N_{\delta_0}[\Sigma(t)h],\mu
\rangle =\delta_0\int_D
N_{\delta_0}[\sigma(t,\cdot)h](x) \mu(dx).
\end{equation}

We can now introduce the limiting equation. It is the one-dimensional
stochastic differential equation
%
%
\begin{equation}
\label{eq2}\quad
\cases{
dv(t)=\hat{F}(t,v(t)) \,dt+\langle\hat{G}(t,v(t)),dw^Q(t)
\rangle_H+\langle\hat{\Sigma}(t),dw^B(t)\rangle_Z,\vspace*{2pt}\cr
v(0)=\langle u_0,\mu\rangle.}
\end{equation}
As the mappings $\hat{F}(t,\cdot)\dvtx\reals\to\reals$ and $\hat
{G}(t,\cdot)\dvtx\reals\to H$ are both Lipschitz continuous, uniformly
with respect to $t \in[0,T]$, for any $T>0$, equation (\ref{eq2})
admits a unique strong solution $v \in L^p(\Omega;C([0,T];\reals))$
for any $p\geq1$ and $T>0$, that is, there exists a unique adapted
process in $L^p(\Omega;C([0,T];\reals))$ which is adapted to the
filtration $\{\mathcal{F}_t\}_{t\geq0}$ such that
\[
v(t)=\langle u_0,\mu\rangle +\int_0^t \hat{F}(s,v(s)) \,ds+\hat
{w}_{A,Q}(v)(t)+\hat{w}_{A,B}(t),
\]
where
\[
\hat{w}_{A,Q}(v)(t):=\int_0^t\langle\hat{G}(s,v(s)),dw^Q(s)\rangle_H,\qquad
\hat{w}_{A,B}(t):=\int_0^t\langle\hat{\Sigma}(s),dw^B(s)\rangle_Z.
\]

Notice that both $\hat{w}_{A,Q}(v)(t)$ and $\hat{w}_{A,B}(t)$ are
$\mathcal{F}_t$-martingales having zero mean.
Moreover, we have
%
%
\begin{equation}
\label{br1}
\E|\hat{w}_{A,Q}(v)(t)|^2=\int_0^t\E|Q \hat{G}(s,v(s))|_H^2 \,ds
\end{equation}
and
%
%
\begin{equation}
\label{br2}
\E|\hat{w}_{A,B}(t)|^2=\int_0^t\E|B \hat{\Sigma}(s)|_Z^2 \,ds.
\end{equation}
In particular, as $w^Q$ and $w^B$ are independent, we have that $\hat
{w}_{A,Q}(v)(t)+\hat{w}_{A,B}(t)$ is an $\mathcal{F}_t$-martingale
having zero mean and covariance
%
%
\begin{equation}
\label{br3}
\int_0^t\bigl(\E|Q \hat{G}(s,v(s))|_H^2+|B \hat{\Sigma
}(s)|_Z^2\bigr) \,ds
\end{equation}
so that there exists some Brownian motion $\beta_t$ defined on some
stochastic basis $(\hat{\Omega}, \hat{\mathcal{F}}, {\hat
\mathcal{F}}_t,\hat{\Pro})$ such that the solution of problem
(\ref{eq2}) coincides in law with the solution of the problem
\[
\cases{
dv(t)=\hat{F}(t,v(t)) \,dt+\Phi(t,v(t)) \,d\beta_t,\cr
v(0)=\langle u_0,\mu\rangle,}
\]
where
%
%
\begin{equation}
\label{br4}
\Phi(t,v)=\bigl(|Q \hat{G}(t,v)|_H^2+|B \hat{\Sigma}(t)|_Z
\bigr)^{1/2}.
\end{equation}

As shown in Remark \ref{rem21}, in the case where the operator $A$ is
self-adjoint, we have
\[
\mu(dx)=\frac1{|D|} \,dx
\]
so that, due to the definition of $\hat{G}(t,v)$, we get
\[
|Q \hat{G}(t,v)|_H^2=\frac1{|D|^2} |Q g(t,\cdot,v)|_H^2
=\frac1{|D|^2}\int_D|[Q g(t,\cdot,v)](x)|^2 \,dx.
\]
Concerning the boundary term, due to (\ref{adjoint}), we have
\begin{eqnarray*}
|B \hat{\Sigma}(t)|_Z^2&=&\frac{\delta_0^2}{|D|^2}\sum_{k \in
\nat}|\langle N_{\delta_0}[\Sigma(t)Bf_k],1\rangle_H|^2\\
&=& \sum_{k\in\nat}
\frac{\delta_0^2}{|D|^2}|\langle[\Sigma(t)Bf_k],N_{\delta_0}^\star
1\rangle_Z|^2\\
&=&\sum_{k \in\nat}
\frac{1}{|D|^2}|\langle f_k,B\sigma(t,\cdot)\rangle_Z|^2=\frac
1{|D|^2}|B\sigma(t,\cdot)|_Z^2\\
&=&\frac1{|D|^2}\int_{\partial D}|[B \sigma(t,\cdot)](\eta)|^2\,
d\eta.
\end{eqnarray*}
Therefore, in the self-adjoint case, we have
\[
\Phi(t,v)=\frac{1}{|D|}\biggl(\int_D|[Q g(t,\cdot,v)](x)|^2
\,dx+\int_{\partial D}|[B \sigma(t,\cdot)](\eta)|^2 \,d\eta
\biggr)^{1/2}.
\]

Now that we have described the candidate limit equation, we prove that
$u_\varepsilon$ in fact converges to its solution.
\begin{Theorem}
\label{mar8}
Assume Hypotheses \ref{H1}, \ref{H2} and \ref{H3}. Then,
for any $u_0 \in H$, $p\geq1$, $T>0$ and $\theta<1$, and for any
$\delta>0$, we have
%
%
\begin{eqnarray}
\label{mar9}
\E\sup_{t \in[\delta,T]}|u_\varepsilon(t)-v(t)|_{H_\mu}^p&\leq&
c_{T,p,\theta} (\varepsilon+\varepsilon^{{p\theta}/2})
(1+|u_0|_{H_\mu}^p)\nonumber\\[-8pt]\\[-8pt]
&&{}+e^{-{\gamma p \delta}/\varepsilon}|u_0|_{H_\mu}^p,\nonumber
\end{eqnarray}
where $v$ is the solution of the one-dimensional problem (\ref{eq2}).
In particular,
\[
\lim_{\varepsilon\to0} \E\sup_{t \in[\delta,T]}|u_\varepsilon
(t)-v(t)|_{H_\mu}^p=0.
\]
\end{Theorem}
\begin{pf}
We have
\begin{eqnarray*}
u_\varepsilon(t)-v(t)&=&(e^{t A/\varepsilon}u_0-\langle u_0,\mu
\rangle )+\int_0^t \bigl(\hat{F}(s,u_\varepsilon(s))-\hat{F}(s,v(s))
\bigr) \,ds\\
&&{}+\int_0^t \bigl\langle \bigl(\hat{G}(s,u_\varepsilon(s))-\hat
{G}(s,v(s))\bigr),dw^Q(s)\bigr\rangle_H+R_\varepsilon(t),
\end{eqnarray*}
where
%
%
\begin{eqnarray}
\label{repsilon}
R_\varepsilon(t)&:=&\int_0^t e^{(t-s) A/\varepsilon
}F(s,u_\varepsilon(s)) \,ds-\int_0^t \hat{F}(s,u_\varepsilon(s))
\,ds\nonumber\\[-8pt]\\[-8pt]
&&{}+w^\varepsilon_{A,Q}(u_\varepsilon)(t)-\hat{w}_{A,Q}(u_\varepsilon
)(t)+w^{\varepsilon}_{A,B}(t)-\hat{w}_{A,B}(t).
\nonumber
\end{eqnarray}
This yields
%
%
\begin{eqnarray}
\label{f100}
&&|u_\varepsilon(t)-v(t)|_{H_\mu}^p\nonumber\\
&&\qquad\leq c_{T,p}\biggl(|e^{t
A/\varepsilon}u_0-\langle u_0,\mu\rangle|_{H_\mu}^p\nonumber\\
&&\qquad\quad\hspace*{24.4pt}{}+\int_0^t |\hat
{F}(s,u_\varepsilon(s))-\hat{F}(s,v(s))|^p \,ds\\
&&\qquad\quad\hspace*{24.4pt}{}+\biggl|\int_0^t \bigl\langle \bigl(\hat{G}(s,u_\varepsilon(s))-\hat
{G}(s,v(s))\bigr),dw^Q(s)\bigr\rangle _H\biggr|^p\nonumber\\
&&\qquad\quad\hspace*{180.1pt}{}+|R_\varepsilon(t)|_{H_\mu}^p\biggr).\nonumber
\end{eqnarray}
Due to the Lipschitz continuity of $\hat{F}(t,\cdot)\dvtx H_\mu\to\reals
$, for any $0\leq t\leq T$, we have
%
%
\begin{eqnarray}
\label{f1}
&&\E\sup_{s \in[0,t]} \int_0^s |\hat{F}(r,u_\varepsilon
(r))-\hat{F}(r,v(r))|^p \,dr\nonumber\\[-8pt]\\[-8pt]
&&\qquad\leq c_{T,p}\int_0^t\E|u_\varepsilon
(r)-v(r)|_{H_\mu}^p \,dr.\nonumber
\end{eqnarray}
Analogously, due to the Lipschitz continuity
of $\hat{G}(t,\cdot)\dvtx H_\mu\to H$ and the Burkholder--Davis--Gundy
inequality, for any $0\leq t\leq T$, we easily obtain
%
%
\begin{eqnarray}
\label{f2}
&&\E\sup_{s \in[0,t]} \biggl|\int_0^s \langle \hat
{G}(r,u_\varepsilon(r))-\hat{G}(r,v(r)),dw^{Q}(r)\rangle_H\biggr|^p\nonumber\\[-8pt]\\[-8pt]
&&\qquad\leq c_{T,p}\int_0^t\E|u_\varepsilon(r)-v(r)|_{H_\mu}^p \,dr.\nonumber
\end{eqnarray}
Then,
thanks to condition (\ref{ergo}), for any $0\leq t\leq T$,
\begin{eqnarray*}
&&\E|u_\varepsilon(t)-v(t)|_{H_\mu}^p\\
&&\qquad\leq c_{T,p}\biggl(e^{-
{\gamma p t}/\varepsilon}|u_0|_{H_\mu}^p+\E\sup_{t \in
[0,T]}|R_\varepsilon(t)|_{H_{\mu}}^p+\int_0^t\E|u_\varepsilon
(s)-v(s)|_{H_\mu}^p \,ds\biggr)
\end{eqnarray*}
and, by comparison, this yields
%
%
\begin{equation}
\label{f1000}
\int_0^t\E|u_\varepsilon(s)-v(s)|_{H_\mu}^p \,ds\leq c_{T,p}
\Bigl(\varepsilon|u_0|_{H_\mu}^p+\E\sup_{t \in[0,T]}|R_\varepsilon
(t)|_{H_{\mu}}^p\Bigr).
\end{equation}
In view of (\ref{f100}), thanks to (\ref{f1}) and (\ref{f2}),
for any $0<\delta<T$, we obtain
\begin{eqnarray*}
&&\E\sup_{t \in[\delta,T]} |u_\varepsilon(t)-v(t)|_{H_\mu}^p\\
&&\qquad\leq c e^{-{\gamma p \delta}/\varepsilon}|u_0|_{H_\mu}^p+
c_{T,p}\int_0^T\E|u_\varepsilon(s)-v(s)|_{H_\mu}^p \,dt\\
&&\qquad\quad{}+c_{T,p} \E
\sup_{t \in[0,T]}|R_\varepsilon(t)|_{H_{\mu}}^p.
\end{eqnarray*}
Therefore, if we show that, for any $T>0$, $p\geq1$ and $\theta\in(0,1)$,
%
%
\begin{equation}
\label{f3}
\E\sup_{t \in[0,T]}|R_\varepsilon(t)|_{H_\mu}^p\leq c_{T,p,\theta}
\varepsilon^{{p\theta}/2} (1+|u_0|_H^p),
\end{equation}
then we can conclude that (\ref{mar9}) holds.
\end{pf}

Due to (\ref{repsilon}), in order to prove (\ref{f3}) and hence
complete the proof of Theorem~\ref{mar8}, we need the following three lemmas.
\begin{Lemma}
\label{mar7}
Assume Hypotheses \ref{H1}, \ref{H2} and \ref{H3}. Then, for any
$T>0$ and $p\geq1$, and for any $\varepsilon\in(0,1]$, we have
%
%
\begin{eqnarray}
\label{limF}
&&\E\sup_{t \in[0,T]}\biggl|\int_0^t e^{(t-s) A/\varepsilon
}F(s,u_\varepsilon(s)) \,ds-\int_0^t\hat{F}(s,u_\varepsilon(s)) \,ds
\biggr|_{H_\mu}^p\nonumber\\[-8pt]\\[-8pt]
&&\qquad\leq c_{T,p} (1+|u_0|_H^p) \varepsilon^p.\nonumber
\end{eqnarray}
\end{Lemma}
\begin{pf}
Due to Hypothesis \ref{H1}, for any $t \in[0,T]$, we have
\begin{eqnarray*}
&&\bigl|e^{(t-s) A/\varepsilon}F(s,u_\varepsilon(s))-\hat
{F}(s,u_\varepsilon(s))\bigr|_{H_\mu}\\
&&\qquad\leq c e^{-{\gamma
(t-s)}/\varepsilon} |F(s,u_\varepsilon(s))|_{H_\mu}\\
&&\qquad\leq c e^{-{\gamma(t-s)}/\varepsilon} |F(s,u_\varepsilon
(s))|_{H}\\
&&\qquad\leq c_{T} e^{-{\gamma(t-s)}/\varepsilon} \Bigl(1+\sup
_{s\leq T}|u_\varepsilon(s)|_{H}\Bigr).
\end{eqnarray*}
This implies that, for any $t \in[0,T]$,
\begin{eqnarray*}
&&
\biggl|\int_0^t e^{(t-s) A/\varepsilon}F(s,u_\varepsilon(s)) \,ds-\int
_0^t\hat{F}(s,u_\varepsilon(s)) \,ds\biggr|_{H_\mu}^p\\
&&\qquad\leq c_{T,p}
\Bigl(1+\sup_{s\leq T}|u_\varepsilon(s)|^p_{H}\Bigr) \biggl(\int_0^t e^{-
{\gamma s}/\varepsilon} \,ds\biggr)^p
\end{eqnarray*}
so that, thanks to (\ref{epsilon}), for any $\varepsilon\in(0,1]$,
we obtain
\begin{eqnarray*}
&&\E\sup_{t \in[0,T]}\biggl|\int_0^t e^{(t-s) A/\varepsilon
}F(s,u_\varepsilon(s)) \,ds-\int_0^t\hat{F}(s,u_\varepsilon(s)) \,ds
\biggr|_{H_\mu}^p\\
&&\qquad\leq c_{T,p}(1+|u_0|^p_H)\varepsilon^p.
\end{eqnarray*}
\upqed\end{pf}
\begin{Lemma}
\label{mar4}
Assume Hypotheses \ref{H1}, \ref{H2} and \ref{H3}, and fix $T>0$,
$p\geq1$ and $\theta<1$.
Then, there exists some constant $c_{T,p,\theta}>0$ such that for any
$\varepsilon\in(0,1]$,
%
%
\begin{equation}
\label{limG}\qquad
\E\sup_{t \in[0,T]}|w^\varepsilon_{A,Q}(u_\varepsilon)(t)-\hat
{w}_{A,Q}(u_\varepsilon)(t)|_{H_\mu}^p\leq c_{T,p,\theta}
\varepsilon^{{p\theta}/2} (1+|u_0|_H^p).
\end{equation}
\end{Lemma}
\begin{pf}
As in the proofs of Lemmas \ref{lemma31} and \ref{lemma33}, we use a
factorization argument. Since $e^{tA} 1=1$, for any $t\geq0$ and $\alpha
>0$, we have
\begin{eqnarray*}
&&w^\varepsilon_{A,Q}(u_\varepsilon)(t)-\int_0^t\langle\hat
{G}(s,u_\varepsilon(s)),dw^Q(s)\rangle_H\\
&&\qquad=\frac{\sin\pi\alpha}{\pi}\int
_0^t (t-s)^{\alpha-1}e^{(t-s) A/\varepsilon}Y_{\varepsilon,\alpha
}(s) \,ds,
\end{eqnarray*}
where
\[
Y_{\varepsilon,\alpha}(s):=\int_0^s(s-r)^{-\alpha}e^{(s-r)
A/\varepsilon
}\Psi(r,u_\varepsilon(r)) \,dw^Q(r)
\]
and, for any $h_1, h_2 \in H$,
\[
\Psi(r,h_1)h_2:=G(r,h_1)h_2-\langle\hat{G}(r,h_1),h_2\rangle_H.
\]
Hence, due to (\ref{mar10}), $e^{tA}$ is a contraction in $H_\mu$ for
any $t\geq0$, and by proceeding as in the proofs of Lemmas \ref
{lemma31} and \ref{lemma33}, for $\alpha<1/p$, we obtain
\begin{eqnarray*}
&&\E\sup_{t \in[0,T]}\biggl|w^\varepsilon_{A,Q}(u_\varepsilon
)(t)-\int_0^t\langle\hat{G}(s,u_\varepsilon(s)),dw^Q(s)\rangle_H
\biggr|^p_{H_\mu}\\
&&\qquad\leq c_{T,p,\alpha} \E\int_0^T\biggl(\int_0^s(s-r)^{-2\alpha}\sum_{k
\in\nat}\lambda_k^2 \bigl|e^{(s-r) A/\varepsilon}\Psi
(r,u_\varepsilon(r))e_k\bigr|_{H_\mu}^2 \,dr\biggr)^{p/2}\,ds.
\end{eqnarray*}
Due to the invariance of $\mu$ and condition (\ref{ergo}), we have
\begin{eqnarray*}
&&\bigl|e^{(s-r) A/\varepsilon}\Psi(r,u_\varepsilon(r))e_k
\bigr|_{H_\mu}\\
&&\qquad=\bigl|e^{(s-r) A/\varepsilon}[G(r,u_\varepsilon
(r))e_k]-\langle \hat{G}(r,u_\varepsilon(r)),e_k\rangle_H\bigr|_{H_\mu}\\
&&\qquad=\bigl|e^{({s-r})/2 A/\varepsilon}\bigl(e^{({s-r})/2
A/\varepsilon}[G(r,u_\varepsilon(r))e_k]\bigr)\\
&&\hspace*{61.6pt}{}-\bigl\langle e^{
({s-r})/2 A/\varepsilon}[G(r,u_\varepsilon(r))e_k],\mu\bigr\rangle
\bigr|_{H_\mu}\\
&&\qquad\leq c e^{-{\gamma(s-r)}/({2\varepsilon})}\bigl|e^{
({s-r})/2 A/\varepsilon}[G(r,u_\varepsilon(r))e_k]\bigr|_{H_\mu}
\end{eqnarray*}
so that
%
%
\begin{eqnarray}
\label{mar5}\hspace*{18pt}
&&\E\sup_{t \in[0,T]}\biggl|w^\varepsilon_{A,Q}(u_\varepsilon
)(t)-\int_0^t\langle \hat{G}(s,u_\varepsilon(s)),dw^Q(s)\rangle_H
\biggr|^p_{H_\mu}\nonumber\\
&&\qquad\leq c_{T,p,\alpha} \E\int_0^T \biggl(\int_0^s(s-r)^{-2\alpha
}e^{-{\gamma(s-r)}/{\varepsilon}}\\
&&\qquad\quad\hspace*{70.5pt}{}\times\sum_{k \in\nat}\lambda_k^2
\bigl|e^{({s-r})/2 A/\varepsilon}[G(r,u_\varepsilon(r))e_k
]\bigr|_{H_\mu}^2 \,dr\biggr)^{p/2} \,ds.\nonumber
\end{eqnarray}
Using the same arguments that were used in the proof of Lemma \ref
{lemma33} [see (\ref{33}) and (\ref{33bis})], for any $0\leq r\leq
s\leq T$, we get
\[
\sum_{k \in\nat}\lambda_k^2 \bigl|e^{({s-r})/2
A/\varepsilon}[G(r,u_\varepsilon(r))e_k]\bigr|_{H_\mu}^2\leq c_T
\biggl[\biggl(\frac\varepsilon{s-r}\biggr)^{ d/({2\zeta})}+ 1\biggr]
\bigl(1+|u_\varepsilon(r)|_H^{2}\bigr),
\]
with $\zeta=\rho/(\rho-2)$ if $d>1$ and with $\zeta=1$ if $d=1$.
Thanks to (\ref{epsilon}), this yields
\begin{eqnarray*}
&&\E\sup_{t \in[0,T]}\biggl|w^\varepsilon_{A,Q}(u_\varepsilon
)(t)-\int_0^t\langle \hat{G}(s,u_\varepsilon(s)),dw^Q(s)\rangle_H
\biggr|^p_{H_\mu}\\
&&\qquad\leq c_{T,p,\alpha}(1+|u_0|_H^p)\varepsilon^{-\alpha p}\biggl(\int
_0^T\biggl[\biggl(\frac\varepsilon t\biggr)^{(2\alpha+{d}/({2\zeta}))}+
1\biggr]e^{-{\gamma t}/{\varepsilon}} \,dt\biggr)^{p/2}.
\end{eqnarray*}
Now, according to the first condition in Hypothesis \ref{H3},
we have $d/2\zeta<1$ so that, for any $\theta<1$,
we can fix $\bar{\alpha}>0$ such that
\[
1-2\bar{\alpha}>\theta,\qquad 2\bar{\alpha}+\frac d{2\zeta}<1.
\]
Then, with a change of variable, we easily obtain
\[
E\sup_{t \in[0,T]}\biggl|w^\varepsilon_{A,Q}(u_\varepsilon)(t)-\int
_0^t\langle \hat{G}(s,u_\varepsilon(s)),dw^Q(s)\rangle_H\biggr|^p_{H_\mu}\leq
c_{T,p,\theta} \varepsilon^{{p\theta}/2}(1+|u_0|_H^p)
\]
for any $p>\bar{p}:=1/\bar{\alpha}$. By the H\"older inequality, we
obtain an analogous estimate for any $p\geq1$ and (\ref{limG}) then follows.
\end{pf}
\begin{Lemma}
\label{mar6}
Assume Hypotheses \ref{H1}, \ref{H2} and \ref{H3}, and fix any
$T>0$, $p\geq1$ and $\theta<1$.
Then, there exists some constant $c_{T,p,\theta}>0$ such that for any
$\varepsilon\in(0,1]$,
%
%
\begin{equation}
\label{limSigma}
\E\sup_{t \in[0,T]}|w^\varepsilon_{A,B}(t)-\hat{w}_{A,B}(t)
|_{H_\mu}^p
\leq c_{T,p,\theta} \varepsilon^{{ p\theta}/2}.
\end{equation}
\end{Lemma}
\begin{pf}
Notice that $(\delta_0-A)e^{tA}1=\delta_0$ for any $t\geq0$. Then, as in
Lemma \ref{lemma31}, by factorization, we obtain
\[
w^\varepsilon_{A,B}(t)-\int_0^t\langle \hat{\Sigma}(s),dw^B(s)
\rangle_Z=\frac{\sin\pi\alpha}{\pi}\int_0^t (t-s)^{\alpha-1}e^{(t-s)
A/\varepsilon}Y_{\varepsilon,\alpha}(s) \,ds,
\]
where
\[
Y_{\varepsilon,\alpha}(s):=\int_0^s(s-r)^{-\alpha}(\delta
_0-A)e^{(s-r)
A/\varepsilon}\Psi(r) \,dw^B(r),
\]
and for any $h \in Z$,
\[
\Psi(r)h:=N_{\delta_0}[\Sigma(r)h]-\frac1{\delta_0}\langle \hat{\Sigma
}(r),h\rangle_Z.
\]
Hence, according to (\ref{mar10}), by arguing as in the proofs of
Lemmas \ref{lemma31} and \ref{lemma33}, for any $p>1/\alpha$, we obtain
\begin{eqnarray*}
&&
\E\sup_{t \in[0,T]} \biggl|w^\varepsilon_{A,B}(t)-\int_0^t
\langle \hat{\Sigma}(s),dw^B(s)\rangle_Z\biggr|_{H_\mu}^p\\
&&\qquad
\leq c_{T,p,\alpha}\int_0^T\biggl(\int_0^s(s-r)^{-2\alpha}\\
&&\qquad\quad\hspace*{63.1pt}{}\times\sum_{k \in
\nat}\theta_k^2 \bigl|(\delta_0-A)e^{(s-r) A/\varepsilon}[\Psi
(r)f_k]\bigr|_{H_\mu}^2 \,dr\biggr)^{p/2} \,ds.
\end{eqnarray*}
Due to the invariance of $\mu$ and to condition (\ref{ergo}), we have
\begin{eqnarray*}
&&
\bigl|(\delta_0-A)e^{(s-r) A/\varepsilon}[\Psi(r)f_k]\bigr|_{H_\mu
}\\
&&\qquad=\bigl|(\delta_0-A)e^{(s-r) A/\varepsilon}N_{\delta_0}[\Sigma
(r)f_k]-\delta
_0\langle N_{\delta_0}[\Sigma(r)f_k],\mu\rangle\bigr|_{H_\mu}\\
&&\qquad=\bigl|e^{({s-r})/{2} A/\varepsilon}\bigl((\delta_0-A)e^{
({s-r})/2 A/\varepsilon}N_{\delta_0}[\Sigma(r)f_k]\bigr)\\
&&\qquad\quad\hspace*{27.37pt}{} - \bigl\langle (\delta
_0-A)e^{({s-r})/2 A/\varepsilon}N_{\delta_0}[\Sigma(r)f_k],\mu
\bigr\rangle\bigr|_{H_\mu}\\
&&\qquad\leq c e^{-{\gamma(s-r)}/({2\varepsilon})}\bigl|(\delta
_0-A)e^{({s-r})/2 A/\varepsilon}N_{\delta_0}[\Sigma(r)f_k]
\bigr|_{H_\mu}.
\end{eqnarray*}
This implies that
\begin{eqnarray*}
&&\E\sup_{t \in[0,T]} \biggl|w^\varepsilon_{A,B}(t)-\int_0^t
\langle \hat{\Sigma}(s),dw^B(s)\rangle_Z\biggr|_{H_\mu}^p\\
&&\qquad\leq c_{T,p,\alpha}\int_0^T\biggl(\int_0^s(s-r)^{-2\alpha}e^{-
{\gamma(s-r)}/{\varepsilon}}\\
&&\qquad\quad\hspace*{64.3pt}{}\times\sum_{k \in\nat}\theta_k^2 \bigl|(\delta
_0-A)e^{({s-r})/2 A/\varepsilon}N_{\delta_0}[\Sigma(r)f_k]
\bigr|_{H_\mu}^2 \,dr\biggr)^{p/2} \,ds
\end{eqnarray*}
and, hence, by proceeding as in the proof of Lemma \ref{lemma31}, we
conclude that
\begin{eqnarray*}
&&\E\sup_{t \in[0,T]} \biggl|w^\varepsilon_{A,B}(t)-\int_0^t
\langle \hat{\Sigma}(s),dw^B(s)\rangle_Z\biggr|_{H_\mu}^p\\
&&\qquad\leq c_{T,p,\alpha,\rho} \varepsilon^{-\alpha p}\biggl(\int_0^T\biggl[
\biggl(\frac\varepsilon s\biggr)^{2\alpha+({d \operatorname{sign} (d-1)+\zeta
})/({2\zeta})+{\rho}/2}+ 1\biggr]e^{-{\gamma s}/{\varepsilon}}\,
ds\biggr)^{p/2},
\end{eqnarray*}
where $\rho$ is a positive constant to be chosen and where $\zeta
=\beta/(\beta-2)$ if $d>1$ and $\zeta=1$ if $d=1$.
Now, as we are assuming $\beta<2d/(d-1)$ when $d\geq2$, for any
$\theta<1$, we can fix $\bar{\alpha}$ and $\bar{\rho}$ both positive
such that
\[
1-2\bar{\alpha}>\theta,\qquad 2\bar{\alpha}+\frac{d \operatorname
{sign}(d-1)+\zeta}{2\zeta}+\frac{\bar{\rho}}2<1.
\]
Then, with a change of variable, for any $p>\bar{p}=1/\bar{\alpha}$,
\[
\E\sup_{t \in[0,T]} \biggl|w^\varepsilon_{A,B}(t)-\int_0^t\langle \hat
{\Sigma}(s),dw^B(s)\rangle_Z\biggr|_{H_\mu}^p\leq c_{T,p} \varepsilon
^{{p \theta}/2}
\]
and this implies (\ref{limSigma}) for any $p\geq1$.
\end{pf}
\eject
\begin{Remark}
\begin{enumerate}
\item Notice that from (\ref{f1000}), we have
%
%
\begin{equation}
\label{notice}
\E|u_\varepsilon-v|_{L^p(0,T;H_\mu)}^p\leq c_{T,p,\theta}
(\varepsilon^{{p\theta}/2}+\varepsilon)(1+|u_0|_{H_\mu}^p)
\end{equation}
so that
\[
\lim_{\varepsilon\to0}\E|u_\varepsilon-v|_{L^p(0,T;H_\mu)}^p=0.
\]
\item If we take $u_0=\langle u_0,\mu\rangle$, then, for any $p\geq1$,
$T>0$ and $\theta<1$, we have the stronger estimate
%
%
\begin{equation}
\label{mar100}
\E\sup_{t \in[0,T]}|u_\varepsilon(t)-v(t)|_{H_\mu}^p\leq
c_{T,p,\theta} \varepsilon^{{p\theta}/2}(1+|u_0|^p).
\end{equation}
\item From the proofs of Lemmas \ref{mar4} and \ref{mar6}, we
easily see that for any $T>0$ and $p\geq1$,
%
%
\begin{equation}
\label{1649}
\sup_{t \in[0,T]}\E|w^\varepsilon_{A,Q}(u_\varepsilon)(t)-\hat
{w}_{A,Q}(t)|_{H_\mu}^p\leq c_{T,p} \varepsilon^{{p}/2}
(1+|u_0|_H^p)
\end{equation}
and
%
%
\begin{equation}
\label{1650}
\sup_{t \in[0,T]}\E|w^\varepsilon_{A,B}(t)-\hat
{w}_{A,B}(t)|_{H_\mu}^p
\leq c_{T,p} \varepsilon^{{p}/2}.
\end{equation}
Then, for any $T>0$ and $p\geq1$,
\[
\sup_{t \in[0,T]}\E|R_\varepsilon(t)|_{H_\mu}^p\leq c_{T,p}
\varepsilon^{p/2}(1+|u_0|_{H_\mu}^p), \qquad\varepsilon
\in(0,1].
\]
Then, by repeating the arguments used in the proof of Theorem \ref
{mar8}, we have
\[
\sup_{t \in[\delta,T]}\E|u_\varepsilon(t)-v(t)|_{H_\mu}^p\leq
c_{T,p} (\varepsilon+\varepsilon^{{p}/2})(1+|u_0|_{H_\mu
}^p)+e^{-{\gamma p \delta}/\varepsilon}|u_0|_{H_\mu}^p.
\]
Moreover, if $u_0=\langle u_0,\mu\rangle$, as in (\ref{mar100}), we have
%
%
\begin{equation}
\label{1651}
\sup_{t \in[0,T]}\E|u_\varepsilon(t)-v(t)|_{H_\mu}^p\leq
c_{T,p} \varepsilon^{{p}/2}(1+|u_0|^p).
\end{equation}
\end{enumerate}
\end{Remark}

\section{Fluctuations around the averaged motion}

In this section, we analyze the fluctuations of the motion
$u_\varepsilon
$ around the averaged motion $v$. More precisely, we will study the
limiting behavior of the random field
%
%
\begin{equation}
\label{fluc}
z_\varepsilon(t,x):=\frac{u_\varepsilon(t,x)-v(t)}{\sqrt{\varepsilon}},\qquad
t \geq0, x \in D,
\end{equation}
as the parameter $\varepsilon$ goes to zero.

In what follows, in addition to Hypothesis \ref{H2}, we shall assume
that the coefficients $f$ and $g$ satisfy the following conditions.
\begin{Hypothesis}
\label{H4}
\begin{enumerate}
\item
The mapping $f(t,x,\cdot)\dvtx\reals\to\reals$ is of class $C^1$, with
Lipschitz continuous derivative, uniformly with respect to $x \in D$
and $t \in[0,T]$, for any $T>0$.
\item The mapping $g$ does not depend on the third variable, that is,
$g(t,x,\eta)=g(t,x)$ for any $t\geq0$, $x \in D$ and $\eta\in
\reals$.
\item For any $x \in D$, the mappings $g(\cdot,x)\dvtx[0,\infty)\to
\reals$ and $\sigma(\cdot,x)\dvtx\reals\to\reals$ are H\"older
continuous of exponent $\alpha>0$ and
%
%
\begin{eqnarray}
\label{rev2}
\sup_{x \in D}[g(\cdot,x)]_{C^\alpha([0,+\infty))}&=&L_g<\infty,\nonumber\\[-8pt]\\[-8pt]
\sup_{\eta
\in\partial D}[\sigma(\cdot,\eta)]_{C^\alpha([0,+\infty
))}&=&L_\sigma<\infty.\nonumber
\end{eqnarray}
\end{enumerate}
\end{Hypothesis}

From Hypothesis \ref{H4}, we easily obtain that the mapping $\hat
{F}(t,\cdot)\dvtx H_\mu\to\reals$ is Fr\'echet differentiable and,
for any $t\geq0$ and $h,k \in H_\mu$, we have
\[
D\hat{F}(t,h)k=\int_D \frac{\partial f}{\partial\xi
}(t,x,h(x))k(x) \mu(dx)=\biggl\langle \frac{\partial f}{\partial\xi
}(t,\cdot,h) k,\mu\biggr\rangle.
\]
Moreover, $D\hat{F}(t,\cdot)\dvtx H_\mu\to H$ is Lipschitz continuous,
uniformly for $t \in[0,T]$.
\begin{Theorem}
\label{teo5.1}
Assume Hypotheses \ref{H1}--\ref{H4}.
Then, for any $t>0$,
%
%
\begin{equation}
\label{1652}
z_\varepsilon(t,x)\rightharpoonup I_0(t,x),\qquad \varepsilon\downarrow0,
\end{equation}
in $H_\mu$, where $I_0(t,x)$ is the Gaussian random field defined for
any $t>0$ and $x \in D$ by
%
%
\begin{eqnarray}
\label{tax7}
I_0(t,x)&:=&\int_0^\infty\Pi e^{sA}G(t)\,dw^Q(s,x)\nonumber\\[-8pt]\\[-8pt]
&&{} + \int_0^\infty\Pi
(\delta_0-A)e^{sA}N_{\delta_0}[\Sigma(t)\,dw^B(s)](x).\nonumber
\end{eqnarray}
[For any $x \in H_\mu$, we have set $\Pi x:=x-\langle x,\mu\rangle$.
Notice that, due to the invariance of $\mu$,
\[
\Pi e^{tA}h=e^{tA} \Pi h,\qquad t\geq0, h \in H_\mu,\qquad \Pi
A h=A \Pi h,\qquad h \in D(A).]
\]
\end{Theorem}

We now define
%
%
\begin{equation}
\label{tax8}
I_G(t):=\int_0^\infty\Pi e^{sA}G(t)\,dw^Q(s)
\end{equation}
and
%
%
\begin{equation}
\label{tax9}
I_\Sigma(t):=\int_0^\infty\Pi(\delta_0-A)e^{sA}N_{\delta_0}[\Sigma
(t)\,dw^B(s)].
\end{equation}
Before proceeding with the proof of Theorem \ref{teo5.1}, it is
important to see that the two terms $I_G(t)$ and $I_\Sigma(t)$
are both well defined in $L^2(\Omega;H_\mu)$ for any $t\geq0$.
\begin{Lemma}
\label{fluc5}
Under Hypotheses \ref{H1}--\ref{H3},
\[
\E|I_G(t)|^2_{H_\mu}<\infty,\qquad t\geq0.
\]
\end{Lemma}
\begin{pf}
Due to the invariance of $\mu$, we have
\[
I_G(t)=\sum_{k=1}^\infty\int_0^\infty\lambda_k
\bigl(e^{sA}[G(t)e_k]-\langle G(t)e_k,\mu\rangle\bigr)\, d\beta_k(s)
\]
so that, by proceeding as in the proof of Lemma \ref{lemma33}, thanks
to (\ref{lambda}),
we have
%
%
\begin{eqnarray}
\label{fluc3}
\E|I_G(t)|^2_{H_\mu}&=&\int_0^\infty\sum_{k=1}^\infty\lambda
_k^2|e^{sA}[G(t)e_k]-\langle G(t)e_k,\mu\rangle|_{H_\mu}^2 \,ds\nonumber\\
&\leq& c \int_0^\infty\Biggl( \sum_{k=1}^\infty\bigl|e^{sA}
\bigl([G(t)e_k]-\langle G(t)e_k,\mu\rangle\bigr)\bigr|_{H_\mu}^2\Biggr)^{1/\zeta
}\\
&&\hspace*{24.2pt}{}\times\sup_{k \in\nat}|e^{sA}[G(t)e_k]-\langle G(t)e_k,\mu
\rangle|_{H_\mu}^{{2(\zeta-1)}/\zeta}|e_k|_\infty^{-
4/\rho}\,ds,\nonumber
\end{eqnarray}
where $\zeta=(\rho-2)/\rho$ and $\rho$ is the constant appearing in
(\ref{lambda}).
Due to (\ref{ergo}) and the invariance of $\mu$, we have
\[
\bigl|e^{sA}\bigl([G(t)e_k]-\langle G(t)e_k,\mu\rangle\bigr)\bigr|_{H_\mu}^2\leq
e^{-\gamma s}|e^{s/ 2 A} [G(t)e_k]|_{H_\mu}^2
\]
so that,
according to (\ref{33}), we have
\[
\Biggl( \sum_{k=1}^\infty\bigl|e^{sA}\bigl([G(t)e_k]-\langle G(t)e_k,\mu
\rangle\bigr)\bigr|_{H_\mu}^2\Biggr)^{1/\zeta}\leq c_t e^{-{\gamma
s}/\zeta} \bigl(s^{- d/({2\zeta})}+1\bigr).
\]
Analogously, according to (\ref{33bis}), we have
\[
\sup_{k \in\nat}|e^{sA}[G(t)e_k]-\langle G(t)e_k,\mu\rangle
|_{H_\mu}^{{2(\zeta-1)}/\zeta}|e_k|_\infty^{-4/\rho}\leq
c_t e^{-{\gamma(\zeta-1)s}/\zeta}
\]
and hence, in view of (\ref{fluc3}), we conclude that
\[
\E|I_G(t)|^2_{H_\mu}\leq c_t \int_0^\infty e^{-\gamma s}
\bigl(s^{- d/({2\zeta})}+1\bigr) \,ds\leq c_t.
\]
\upqed\end{pf}

As far as $I_\Sigma$ is concerned, we have the following, analogous, result.
\begin{Lemma}
\label{fluc6}
Under Hypotheses \ref{H1}--\ref{H3}
\[
\E|I_\Sigma(t)|^2_{H_\mu}<\infty,\qquad t\geq0.
\]
\end{Lemma}
\begin{pf}
Due to the invariance of $\mu$, we have
\[
I_\Sigma(t)=\sum_{k=1}^\infty\int_0^\infty\theta_k \bigl((\delta
_0-A)e^{sA}N_{\delta_0}[\Sigma(t)f_k]-\delta_0\langle N_{\delta_0}[\Sigma
(t)f_k],\mu\rangle\bigr) \,d\hat{\beta}_k(s).
\]
Using the same arguments used in Lemma \ref{mar6}, due to (\ref
{ergo}) and the invariance of~$\mu$,
we have
\begin{eqnarray*}
&&|(\delta_0-A)e^{sA}N_{\delta_0}[\Sigma(t)f_k]-\delta_0\langle N_{\delta
_0}[\Sigma
(t)f_k],\mu\rangle|_{H_\mu}^2\\
&&\qquad\leq c e^{-\gamma s}|(\delta
_0-A)e^{s/2 A}N_{\delta_0}[\Sigma(t)f_k]|_{H_\mu}^2
\end{eqnarray*}
and then, as in the proof of Lemma \ref{lemma31}, due to (\ref
{theta}), we get
\begin{eqnarray*}
\E|I_\Sigma(t)|_{H_\mu}^2&\leq& c \int_0^\infty e^{-\gamma s}
\Biggl( \sum_{k=1}^\infty|(\delta_0-A)e^{ s/2 A}N_{\delta_0}[\Sigma
(t)f_k]|_{H_\mu}^2\Biggr)^{1/\zeta}\\
&&\hspace*{23.8pt}{}\times\sup_{k \in\nat}|(\delta_0-A)e^{ s/2 A}N_{\delta
_0}[\Sigma(t)f_k]|_{H_\mu}^{{2(\zeta-1)}/\zeta} \,ds.
\end{eqnarray*}
By using (\ref{quadrato}) and (\ref{3}), this allows us to conclude that
for some $\bar{\rho}>0$ such that $(d+\zeta)/2\zeta+\bar{\rho}/2<1$,
\[
\E|I_\Sigma(t)|_{H_\mu}^2\leq c_t\int_0^\infty e^{-\gamma s}
\bigl(s^{-({d+\zeta})/({2\zeta})+{\bar{\rho}}/2}+1\bigr)\, ds<+\infty.
\]
\upqed\end{pf}

\subsection[Proof of Theorem 5.1]{Proof of Theorem \protect\ref{teo5.1}}

It is immediate to check that for any $t\geq0$,
\[
z_\varepsilon(t)=\int_0^t D\hat{F}(s,v(s))z_\varepsilon(s)
\,ds+R_\varepsilon(t)+I_\varepsilon(t),
\]
where
\begin{eqnarray*}
R_\varepsilon(t):\!&=&\frac1{\sqrt{\varepsilon}}(e^{
t/{\varepsilon} A} u_0-\langle u_0,\mu\rangle)\\
&&{}+\frac1{\sqrt{\varepsilon
}}\int_0^t\bigl(e^{(t-s) A/\varepsilon}F(s,u_\varepsilon(s))-\hat
{F}(s,u_\varepsilon(s))\bigr)\,ds\\
&&{}+\int_0^t
\int_0^1\bigl[D\hat{F}\bigl(s,v(s)+\theta\bigl(u_\varepsilon(s)-v(s)\bigr)\bigr)-D\hat
{F}(s,v(s))\bigr]z_\varepsilon(s) \,ds \,d\theta\\
&=&\!:\sum_{i=1}^3
R_{\varepsilon,i}(t)
\end{eqnarray*}
and
%
%
\begin{equation}
\label{tax5}
I_\varepsilon(t):=\frac1{\sqrt{\varepsilon}}\bigl(w^\varepsilon
_{A,Q}(t)-\hat{w}_{A,Q}(t)\bigr)+\frac1{\sqrt{\varepsilon}}
\bigl(w^\varepsilon_{A,B}(t)-\hat{w}_{A,B}(t)\bigr).
\end{equation}

Due to (\ref{dubbio}), we have
%
%
\begin{equation}
\label{prr1}
|R_{\varepsilon,1}(t)|_{H_\mu}\leq\frac c{\sqrt{\varepsilon}}
e^{-{\gamma t}/{\varepsilon}} |u_0|_{H_\mu}.
\end{equation}
For $R_{\varepsilon,2}(t)$,
with a change of variables, due to (\ref{ergo}), we have, for any $t
\in[0,T]$,
\begin{eqnarray*}
|R_{\varepsilon,2}(t)|_{H_\mu}&\leq&\frac{c}{\sqrt{\varepsilon}}\int
_0^{t}e^{-{\gamma(t-s)}/\varepsilon}|F(s,u_\varepsilon(s))
|_{H_\mu} \,ds\\
&\leq&\frac{c_t}{\sqrt{\varepsilon}}\Bigl(1+\sup_{s \in
[0,t]}|u_\varepsilon(s)|_{H_\mu}\Bigr)\int_0^{t}e^{-{\gamma
s}/\varepsilon} \,ds\\
&\leq& c_t {\sqrt{\varepsilon}}\Bigl(1+\sup_{s \in
[0,t]}|u_\varepsilon(s)|_{H_\mu}\Bigr)
\end{eqnarray*}
and then, thanks to (\ref{epsilon}), we get
%
%
\begin{equation}
\label{i1}
\E\sup_{t \in[0,T]}|R_{\varepsilon,2}(t)|_{H_\mu}\leq c_T \sqrt
{\varepsilon} (1+|u_0|_{H_\mu}).
\end{equation}
Finally, for $R_{\varepsilon,3}(t)$,
due to the Lipschitz continuity of $D\hat{F}(s,\cdot)\dvtx H_\mu\to H$,
uniform with respect to $s \in[0,t]$, and estimate
(\ref{notice}) with $p=2$ and $\theta\in(1/2,1)$, we get
%
%
\begin{eqnarray}
\label{fluc4}
\E|R_{\varepsilon,3}(t)|_{H_\mu}&\leq&\frac{c_t}{\sqrt{\varepsilon
}} \int_0^T\E|u_\varepsilon(s)-v(s)|_{H_\mu}^2 \,ds\nonumber\\[-8pt]\\[-8pt]
&\leq& c_T
\varepsilon^{\theta-1/2} (1+|u_0|^2).\nonumber
\end{eqnarray}
Therefore, collecting together (\ref{prr1}), (\ref{i1}) and (\ref
{fluc4}), we can conclude that for any $T>0$ and $\varepsilon\in(0,1]$,
%
%
\begin{equation}
\label{tax4}
\E|R_\varepsilon(t)|_{H_\mu}\leq\frac c{\sqrt{\varepsilon}}
e^{-{\gamma t}/{\varepsilon}} |u_0|_{H_\mu}+c_T
(1+|u_0|^2_{H_\mu}) \varepsilon^{\theta-1/2},\qquad t
\in[0,T].\hspace*{-36pt}
\end{equation}

Next, for any $\varepsilon>0$, we introduce the problem
\[
\zeta(t)=\int_0^t D\hat{F}(s,v(s))\zeta(s) \,ds+I_\varepsilon(t),
\]
where $I_\varepsilon(t)$ is the process introduced in (\ref{tax5}). For
any $\varepsilon>0$, we denote by $\zeta_\varepsilon$ its unique solution.

We have the following result, whose proof is postponed to the end of
this section.
\begin{Lemma}
\label{tax10}
Under Hypotheses \ref{H1}--\ref{H4}, for any $t>0$, we have
\[
\zeta_\varepsilon(t)\rightharpoonup I_0(t),\qquad \varepsilon
\downarrow0,
\]
in $H_\mu$, where $I_0(t)$ is the $H_\mu$-valued Gaussian vector
field defined in (\ref{tax7}).
\end{Lemma}

Now, for any $\varepsilon>0$ and $t>0$, we define $\rho_\varepsilon
(t):=z_\varepsilon(t)-\zeta_\varepsilon(t)$. We have
\[
\rho_\varepsilon(t)=\int_0^t D\hat{F}(s,v(s))\rho_\varepsilon
(s) \,ds+R_\varepsilon(t)
\]
so that
\[
\E|\rho_\varepsilon(t)|_{H_\mu}\leq c_T\int_0^t\E
|\rho_\varepsilon(s)|_{H_\mu}+\E|R_\varepsilon(t)|_{H_\mu}.
\]
By comparison, we get
\[
\E|\rho_\varepsilon(t)|_{H_\mu}\leq c_T \E
|R_{\varepsilon}(t)|_{H_\mu}+c_T\int_0^t \E|R_{\varepsilon
}(s)|_{H_\mu} \,ds
\]
and, thanks to (\ref{tax4}), this implies that
\begin{eqnarray*}
\E|\rho_\varepsilon(t)|_{H_\mu}&\leq&
\frac{c_T}{\sqrt{\varepsilon}} e^{-{\gamma t}/{\varepsilon}}
|u_0|_{H_\mu}+c_T (1+|u_0|^2_{H_\mu}) \varepsilon^{\theta
-1/2}\\
&&{}+
\frac{c_T}{\sqrt{\varepsilon}}\int_0^t e^{-{\gamma s}/{\varepsilon
}} \,ds\,|u_0|_{H_\mu}.
\end{eqnarray*}
Hence, we can conclude that for any $t>0$,
\[
\lim_{\varepsilon\to0}\E|z_\varepsilon(t)-\zeta_\varepsilon(t)
|_{H_\mu}=\lim_{\varepsilon\to0}\E|\rho_\varepsilon(t)
|_{H_\mu}=0
\]
so that, in view of Lemma \ref{tax10}, Theorem \ref{teo5.1} is proved.

\subsubsection[Proof of Lemma 5.4]{Proof of Lemma \protect\ref{tax10}}

For any $x \in D$ and $t>0$, we have
\[
\zeta_\varepsilon(t,x)=\int_0^t \int_D\frac{\partial f}{\partial\xi
}(s,y,v(s))\zeta_\varepsilon(s,y) \mu(dy) \,ds+I_\varepsilon(t,x).
\]
Then, if we multiply both sides above by $\partial f/\partial\xi
(t,x,v(t))$ and integrate in $x$ with respect to the measure
$\mu$, we get
\[
\Psi_\varepsilon(t)=H(t)\int_0^t\Psi_\varepsilon(s) \,ds+K_\varepsilon(t),
\]
where
\[
\Psi_\varepsilon(t):=\int_D \frac{\partial f}{\partial\xi
}(t,x,v(t))\zeta_\varepsilon(t,x) \mu(dx)
\]
and
\begin{eqnarray*}
H(t)&:=&\int_D \frac{\partial f}{\partial\xi}(t,x,v(t)) \mu(dx),\\
K_\varepsilon(t)&:=&\int_D\frac{\partial f}{\partial\xi
}(t,x,v(t))I_\varepsilon(t,x) \mu(dx).
\end{eqnarray*}
It is then immediate to check that
\[
\int_0^t\Psi_\varepsilon(s) \,ds=\int_0^t\exp\biggl(\int_s^t H(r)
\,dr\biggr)K_\varepsilon(s) \,ds
\]
so that
\[
\zeta_\varepsilon(t,x)=\int_0^t H(t,s)K_\varepsilon(s) \,ds+I_\varepsilon(t,x),
\]
where
\[
H(t,s):=\exp\biggl(\int_s^t H(r) \,dr\biggr).
\]

\textit{Step} 1. We show that for any $t\geq0$,
%
%
\begin{equation}
\label{tax11}
\lim_{\varepsilon\to0}\E\biggl|\int_0^tH(t,s)K_\varepsilon(s) \,ds\biggr|^2=0.
\end{equation}

Due to (\ref{tax5}) and the stochastic Fubini theorem, we have
\begin{eqnarray*}
&&\int_0^tH(t,s)K_\varepsilon(s) \,ds\\
&&\qquad=\frac1{\sqrt{\varepsilon}}\sum_{k=0}^\infty\int_0^t\int_\sigma
^tH(t,s)\biggl\langle \frac{\partial f}{\partial\xi}(s,\cdot,v(s)),\\
&&\qquad\quad\hspace*{102.4pt}e^{(s-\sigma
)A/\varepsilon}\Pi[G(\sigma)Qe_k]\biggr\rangle_{H_\mu}\,ds \,d\beta
_k(\sigma)\\
&&\qquad\quad{}+\frac1{\sqrt{\varepsilon}}\sum_{k=0}^\infty\int_0^t\int_\sigma
^tH(t,s)\biggl\langle \frac{\partial f}{\partial\xi}(s,\cdot,v(s)),\\
&&\qquad\quad\hspace*{115.1pt}(\delta
_0-A)e^{(s-\sigma)A/\varepsilon}\\
&&\qquad\quad\hspace*{115.1pt}{}\times\Pi[N_{\delta_0}(\Sigma(\sigma
)Bf_k)]\biggr\rangle_{H_\mu} \,ds \,d\hat{\beta}_k(\sigma).
\end{eqnarray*}
Then, as $w^Q$ and $w^B$ are independent and $\partial f/\partial\xi$
is uniformly bounded, we get
\begin{eqnarray*}
&&\E\biggl|\int_0^tH(t,s)K_\varepsilon(s) \,ds\biggr|^2\\
&&\qquad\leq\frac{\kappa_t}\varepsilon\int_0^t \sum_{k=0}^\infty
\biggl(\int_\sigma^t e^{\kappa_t(t-s)}
\bigl|e^{(s-\sigma) A/\varepsilon}\Pi[G(\sigma)Qe_k]\bigr|_{H_\mu
} \,ds\biggr)^2 \,d\sigma\\
&&\qquad\quad{}+\frac{\kappa_t}\varepsilon\int_0^t \sum_{k=0}^\infty\biggl(\int
_\sigma^t e^{\kappa_t(t-s)}
\bigl|(\delta_0-A)e^{(s-\sigma) A/\varepsilon}\\
&&\qquad\quad\hspace*{110.7pt}{}\times\Pi[N_{\delta_0}
(\Sigma(\sigma)Bf_k)]\bigr|_{H_\mu} \,ds\biggr)^2 \,d\sigma\\
&&\qquad=:\frac{\kappa_t}\varepsilon\int_0^t\bigl(J_{\varepsilon,1}(t,\sigma
)+J_{\varepsilon,2}(t,\sigma)\bigr) \,d\sigma.
\end{eqnarray*}

For the first term $J_{\varepsilon,1}$, in view of (\ref{ergo}), for
any $\alpha\in(0,2)$, we have
\begin{eqnarray*}
&&
J_{\varepsilon,1}(t,\sigma)\\
&&\qquad\leq\sum_{k=0}^\infty\lambda_k^2 \biggl(\int
_\sigma^t e^{\kappa_t(t-s)}e^{-{\gamma(s-\sigma)}/({2\varepsilon})}
\bigl|e^{(s-\sigma) A/({2\varepsilon})}\Pi[G(\sigma)e_k]\bigr|_{H_\mu
} \,ds\biggr)^2\\
&&\qquad\leq\biggl(\int_\sigma^t e^{\kappa_t(2-\alpha)(t-s)}e^{-{\gamma
(2-\alpha)(s-\sigma)}/({2\varepsilon})} \,ds\biggr)^{2/({2-\alpha})}\\
&&\qquad\quad{}\times\sum
_{k=0}^\infty\lambda_k^2 \biggl(\int_\sigma^t
\bigl|e^{(s-\sigma) A/({2\varepsilon})}[G(\sigma)e_k]\bigr|_{H_\mu
}^{({2-\alpha})/({1-\alpha})}
\,ds\biggr)^{{2(1-\alpha)}/({2-\alpha})}\\
&&\qquad\leq c_t \varepsilon^{2/({2-\alpha})} \sum_{k=0}^\infty\lambda_k^2
\biggl(\int_\sigma^t
\bigl|e^{(s-\sigma) A/({2\varepsilon})}[G(\sigma)e_k]\bigr|_{H_\mu
}^{({2-\alpha})/({1-\alpha})} \,ds\biggr)^{{2(1-\alpha)}/({2-\alpha})}.
\end{eqnarray*}
Then, if we set $\zeta=\rho/(\rho-2)$, by using the H\"older
inequality for infinite series, we get
\begin{eqnarray*}
&&
J_{\varepsilon,1}(t,\sigma)\\
&&\qquad\leq c_t \varepsilon^{2/({2-\alpha})}
\kappa_Q^{2/\rho}\\
&&\qquad\quad{}\times\Biggl(\int_\sigma^t\Biggl( \sum_{k=0}^\infty
\bigl|e^{(s-\sigma) A/({2\varepsilon})}[G(\sigma)e_k]\bigr|_{H}^{2\zeta
}\\
&&\qquad\quad\qquad\quad\hspace*{65.8pt}{}\times|e_k|_\infty^{-4/({\rho-2})}\Biggr)^{1/\zeta({2-\alpha
})/({2(1-\alpha)})} \,ds\Biggr)^{{2(1-\alpha)}/({2-\alpha})}
\end{eqnarray*}
and, by proceeding as in the proof of Lemma \ref{lemma33}, we conclude
that for $\varepsilon\in(0,1]$,
\[
J_{\varepsilon,1}(t,\sigma)\leq c_t \varepsilon^{2/({2-\alpha})}
\kappa_Q^{2/\rho}
\biggl(\int_\sigma^t\bigl( (s-\sigma)^{-d/({2\zeta})({2-\alpha})/({2(1-\alpha
)})}+1\bigr) \,ds\biggr)^{{2(1-\alpha)}/({2-\alpha})}.
\]
Now, in view of Hypothesis \ref{H3}, we have $d/2\zeta<1$ and can fix
$\bar{\alpha}_1>0$ such that
\[
\frac d{2\zeta}\frac{2-\bar{\alpha}_1}{2(1-\bar{\alpha}_1)}<1
\]
and then
%
%
\begin{equation}
\label{tax20}
\frac{\kappa_t}{\varepsilon}\int_0^t J_{\varepsilon,1}(t,\sigma) \,d\sigma
\leq c_{t,\bar{\alpha}_1} \varepsilon^{{\bar{\alpha}_1}/({2-\bar
{\alpha
}_1})},\qquad \varepsilon\in(0,1], t\geq0.
\end{equation}
The same arguments can be repeated for the term $J_{\varepsilon,2}$, so
we can find some $\bar{\alpha}_2>0$ such that
\[
\frac{\kappa_t}{\varepsilon}\int_0^t J_{\varepsilon,2}(t,\sigma) \,d\sigma
\leq c_{t,\bar{\alpha}_2} \varepsilon^{{\bar{\alpha}_2}/({2-\bar
{\alpha
}_2})},\qquad \varepsilon\in(0,1], t\geq0.
\]
This, together with (\ref{tax20}), implies that
\[
\E\biggl|\int_0^tH(t,s)K_\varepsilon(s) \,ds\biggr|^2\leq c_t \varepsilon
^{\gamma},\qquad \varepsilon\in(0,1], t\geq0,
\]
where
\[
\gamma=\frac{\bar{\alpha}_1\wedge\bar{\alpha}_2}{2-\bar{\alpha}_1\wedge
\bar{\alpha}_2},
\]
so (\ref{tax11}) follows.

\textit{Step} 2. We show that for any fixed $t>0$,
%
%
\begin{equation}
\label{tax21}
I_\varepsilon(t)\rightharpoonup I_0(t), \qquad\varepsilon\downarrow0.
\end{equation}

With a change of variable, we have
\begin{eqnarray*}
I_\varepsilon(t)&=&\frac1{\sqrt{\varepsilon}}\biggl( \int_0^t
e^{(t-s)A/\varepsilon}\Pi[G(s)\,dw^Q(s)]\\
&&\hspace*{22.4pt}{}+\int_0^t
(\delta_0-A)e^{(t-s) A/\varepsilon}\Pi[N_{\delta_0}(\Sigma
(s)\,dw^B(s))]\biggr)\\
&=&\int_0^{t/\varepsilon} e^{r A}\Pi[G(t-\varepsilon
r)\,dw_{\varepsilon,t}^Q(r)]\\
&&{}+\int_0^{t/\varepsilon}
(\delta_0-A)e^{r A}\Pi\bigl[N_{\delta_0}\bigl(\Sigma(t-\varepsilon
r)\,dw_{\varepsilon,t}^B(r)\bigr)\bigr],
\end{eqnarray*}
where
\[
w^Q_{\varepsilon,t}(r)=\frac1{\sqrt{\varepsilon}}
\bigl(w^Q(t)-w^Q(t-\varepsilon r)\bigr),\qquad w^B_{\varepsilon,t}(r)=\frac
1{\sqrt{\varepsilon}}\bigl(w^B(t)-w^B(t-\varepsilon r)\bigr).
\]
This means that for any $\varepsilon>0$ and $t>0$,
\[
\mathcal{L}(I_\varepsilon(t))=\mathcal{L}(\hat{I}_\varepsilon(t)),
\]
where
\begin{eqnarray*}
\hat{I}_\varepsilon(t)&:=&\int_0^{t/\varepsilon} e^{r A}\Pi
[G(t-\varepsilon r)\,dw^Q(r)]\\
&&{}+\int_0^{t/\varepsilon}
(\delta_0-A)e^{r A}\Pi\bigl[N_{\delta_0}\bigl(\Sigma(t-\varepsilon r)\,dw^B(r)
\bigr)\bigr].
\end{eqnarray*}
Thus, in order to obtain (\ref{tax21}), it is sufficient to prove
%
%
\begin{equation}
\label{tax22}
\lim_{\varepsilon\to0}\E|\hat{I}_\varepsilon(t)-I_0(t)
|_{H_{\mu}}^2=0.
\end{equation}

We have
\begin{eqnarray*}
\hat{I}_\varepsilon(t)-I_0(t)&=&\int_0^{t/\varepsilon}e^{r
A}\Pi\bigl[\bigl(G(t-\varepsilon r)-G(t)\bigr)\,dw^Q(r)\bigr]\\
&&{}+\int_0^{t/\varepsilon}
(\delta_0-A)e^{r A}\Pi\bigl[N_{\delta_0}\bigl(\bigl(\Sigma(t-\varepsilon
r)-\Sigma(t)\bigr)\,dw^B(r)\bigr)\bigr]\\
&&{}-\int_{t/\varepsilon}^\infty e^{r A}\Pi[G(t)\,dw^Q(r)
]\\
&&{}-\int_{t/\varepsilon}^\infty
(\delta_0-A)e^{r A}\Pi[N_{\delta_0}(\Sigma(t)\,dw^B(r))]\\
&=&\!:\sum
_{i=1}^4 J_{\varepsilon,i}(t).
\end{eqnarray*}

With the same arguments used several times throughout the paper, we have
\[
\E|J_{\varepsilon,1}(t)|_{H_\mu}^2\leq c \int_0^{t/\varepsilon
}e^{-\gamma s}\bigl(s^{-{d}/({2\zeta})}+1\bigr) |g(t-\varepsilon
s,\cdot)-g(t,\cdot)|_{H_\mu}^2 \,ds.
\]
Then, due to Hypothesis \ref{H4}, we have
%
%
\begin{equation}
\label{fluc9}
\E|J_{\varepsilon,1}(t)|_{H_\mu}^2\leq c_t \varepsilon^{2\alpha} \int
_0^{\infty}e^{-\gamma s}\bigl(s^{-{d}/({2\zeta})}+1\bigr)s^{2\alpha}
\,ds\leq c_t \varepsilon^{2\alpha}.
\end{equation}
Analogously, we have
%
%
\begin{equation}
\label{fluc92}\E|J_{\varepsilon,2}(t)|_{H_\mu}^2\leq c_t \varepsilon
^{2\alpha}.
\end{equation}
Concerning $J_{\varepsilon,3}(t)$, we have
\begin{eqnarray*}
\E|J_{\varepsilon,3}(t)|_{H_\mu}^2 &\leq& c\int_{t/\varepsilon
}^\infty e^{-\gamma s}\bigl(s^{-{d}/({2\zeta})}+1\bigr) \,ds\,
|g(t,\cdot)|_{H_\mu}^2\\
&\leq& c_t\int_{t/\varepsilon}^\infty
e^{-\gamma s}\bigl(s^{-{d}/({2\zeta})}+1\bigr) \,ds
\end{eqnarray*}
so that
%
%
\begin{equation}
\label{fluc10}
\lim_{\varepsilon\to0}\E|J_{\varepsilon,3}(t)|_{H_\mu}^2=0.
\end{equation}
In an identical way, we can show that
\[
\lim_{\varepsilon\to0}\sup_{t \in[\delta,T]}\E|J_{\varepsilon
,4}(t)|_{H_\mu}^2=0
\]
and this, together with
(\ref{fluc9}), (\ref{fluc92}) and (\ref{fluc10}), implies (\ref{tax22}).

\section*{Acknowledgments}
We would like to thank the two anonymous referees
who read the first version of our paper for their interesting and
useful remarks and suggestions.

%

%
\printaddresses

\end{document}